\documentclass{elsart}

\usepackage{amssymb}
\usepackage{graphicx}
\usepackage{enumitem}
\usepackage{subfig}
\usepackage{color}

\newcommand{\be}{\begin{equation}}
\newcommand{\ee}{\end{equation}}
\newcommand{\p}{\partial}

\begin{document}
\begin{frontmatter}



\title{A new parallel solver suited for  
arbitrary semilinear parabolic partial differential equations based on 
generalized random trees} 
\author[label1]{Juan A. Acebr\'{o}n}
\ead{juan.acebron@ist.utl.pt}
and
\author[label1]{\'Angel Rodr\'{\i}guez-Rozas}
\ead{angel.rodriguez@ist.utl.pt}
\address[label1]{Center for Mathematics and its Applications,
Department of Mathematics, Instituto Superior T\'ecnico
Av. Rovisco Pais 1049-001 Lisboa, Portugal}

\begin{abstract}
A probabilistic representation for  
initial value semilinear parabolic problems based on generalized
random trees has been derived. Two different strategies have been proposed,
both requiring generating suitable random trees combined with a Pade
approximant for approximating accurately a given divergent series. Such series are obtained by summing the partial contribution to the solution coming from trees with arbitrary number of branches. The new representation greately expands the class of problems amenable to be solved probabilistically, and was used
successfully to develop a generalized probabilistic domain decomposition 
method. Such a method has been shown to be suited for massively parallel
computers, enjoying full scalability and fault tolerance. Finally, a few numerical examples are given to illustrate the remarkable performance of the algorithm, comparing the results with those obtained with a classical method.

\end{abstract}

\end{frontmatter}

\section{Introduction}\label{introd}

Nowadays, most of the more popular numerical methods developed 
for solving partial differential equations (PDE) are based as a rule in generating a
computational mesh, discretizing the given problem using the nodes of a computational mesh
as discretization nodes, and solving the ensuing linear algebra problem for such
nodes. As an alternative, the so-called meshless methods have started to be exploited more
recently to avoid the need of building a computational mesh, but unfortunately they cannot exclude solving a corresponding linear algebra problem. 
Moreover, both numerical methods share also an important disadvantage, consisting namely in the impossibility of computing the solution of the problem at a single point inside the domain. Clearly this is due to the globally coupled nature
behind these numerical methods. In the past such a feature did not constitute
any serious limitation of the methods, being implemented frequently in sequential computers, but currently with the advent of parallel computers, this can be seen as a 
strong limiting factor to the overall efficiency of the corresponding 
numerical algorithms.

In fact, there are three sources of problem that have been observed. First, the tightly 
coupled nature of the algorithms induces a strong communication among the large number of processors currently present in the more advanced high performance supercomputing, and thus reducing the effective performance. Second, the chance to get a failure in one or several processors during the computation time increases with the number of processors involved. There is indeed a non negligible probability that a small percentage of processors or the network connecting them experience a failure. Most existing algorithms 
simply stop and are aborted as a consequence of such failures, and the proposed remedies usually require some kind of storing and
restarting procedures, which degrades seriously the overall performance. Finally, in case of the grid computing and heterogeneous distributed computing, things can be even worse due to the high degree of systems heterogeneity and high network latency.  Three major and recently published studies about the current source of problems in scientific computing can be found in \cite{Kogge}, \cite{Sarkar}, and \cite{Stevens}.

A successful alternative to such traditional methods consists of developing parallel numerical algorithms
based on probabilistic numerical methods \cite{Acebron1,Acebron2,Angel1,Angel2,Vilela,Arbenz,Ramirez}. Such methods are inherently parallel, and naturally fault-tolerant, hence overcome all the obstacles
mentioned above. Moreover, they allow to compute the solution at a single point inside the domain without the need of a computational mesh for
solving the entire problem. However, they are not recommended to be applied in every point of the mesh due to the slow convergence
rate, but rather to be combined with a classical Domain Decomposition \cite{Chan,Quarteroni} to compute merely the solution in a few points along some interfaces inside the
domain. This is essentially the method called Probabilistic Domain 
Decomposition (PDD) proposed for the first time in \cite{Acebron1,Acebron2} for solving 
linear elliptic problems, and generalized further to deal with some particular semilinear parabolic problems in \cite{Angel1,Angel2}. In short, the idea consists of generating only few interfacial values using probabilistic methods along a given possibly artificial interfaces inside the domain, obtaining approximate values interpolating on such interfaces, and then use such values as boundary data in order to split the original problem into fully decoupled sub-problems.

The aforementioned problems, successfully resolved by the PDD method, are very common and generally recognized by the global high performance computing community. At the present time, the situation becomes even more dramatic because of the number of scientific and engineering challenging problems requiring levels of petaflops, even exaflops, thus computing performance tends to increase very fast. A compelling and comprehensive reference about this fact can be found in the recently published roadmap of the \textit{International Exascale Software Project} (IESP) \cite{Dongarra}. In this project a worldwide joint effort is being conducted trying to overcome the current problems looking for a high quality computational environment for petascale/exascale systems. Incidentally in \cite{Dongarra}, it is claimed the need of rethinking completely new algorithms, and in particular it was suggested to reconsider the use of Monte Carlo based approaches, 
which is actually the case of the PDD method.

A key ingredient for implementing any PDD method requires having a probabilistic
representation of the solution, since it will allow to compute the solution at
the interfacial values.
Probabilistic representations do exist for some elementary {\it semilinear} 
parabolic equations. Indeed, in \cite{McKean} H.P. McKean derived the 
representation formula
\begin{equation}
   u(x,t) = E[\prod_{i=1}^{k_t(\omega)} f(x_i(\omega,t))]     \label{prod-McK}
\end{equation}
for the KPP equation
\begin{equation}
  u_t = u_{xx} + u (u - 1),    \quad   x\in{\bf R},  \  t > 0,
                                                            \label{KPP-1D}
\end{equation}
subject to the initial value $u(x,0) = f(x)$;
see also \cite{Freidlin,Milstein,regnier}. It is understood that in (\ref{prod-McK}) the point
$x_i(\omega,t)$ is the position of the $i$th branch of a branching stochastic process surviving at 
time $t$, $\omega$ denoting the chance variable. The quantity $k_t(\omega)$
is the random number of descendants at time $t$.
  A similar representation has been recently found in \cite{Angel1,Angel2} for the solution of a more general semilinear parabolic problem, given by
\begin{equation}
   \frac{\partial u}{\partial t} = Lu - c u + \sum_{j=2}^{m} \alpha_j u^j,
                                                          \label{nonlinear_prev}
\end{equation}
where $L$ is a general linear elliptic operator,  say $L := 
a_{ij}({\bf x},t) \p_i\p_j + b_i({\bf x},t) \p_i$ (using the summation 
convention), with continuous bounded coefficients, $m \geq 2$ is an integer, $\alpha_j \geq 0$, $\sum_{j=2}^m \alpha_j = 1$,
and $c$ is a positive constant. Such a representation is based on 
generating {\it branching} diffusion processes, associated with the elliptic 
operator in Eq.~(\ref{nonlinear_prev}), and governed by an exponential random time, 
$S$, with probability density $p(S) = c \, \exp(-c S)$. 

In this paper this method is extended to deal with a wider class of semilinear parabolic problems, whose general form now is given by
\begin{eqnarray}
\label{nonlinear_new}
   \frac{\partial u}{\partial t} = Lu + f(u,x,t), \ x \in {\bf R}^n, \ t > 0 \\ \nonumber
	 u(x,0) = g(x),
\end{eqnarray}
where
\begin{equation}
   f(u,x,t) = \sum_{j=2}^m c_j(x,t) u^{j}. \nonumber
\end{equation}

It is worth to observe that this generalizes further the previous representation obtained in \cite{Angel1,Angel2}, since it accounts for
the following aspects: A constant potential term such as $- c u$ is not required anymore; the coefficients multiplying the nonlinear terms,  $c_j(x,t)$, can be now chosen arbitrarily, hence overcoming 
the constraint imposed in the previous representation consisting in $\sum_{j=2}^m c_j(x,t)=1$, and finally the initial data $g(x)$ may now be chosen negative, or greater than $1$.

Moreover, using such a generalized probabilistic representation, the PDD me\-thod
will be generalized further increasing notably the type of semilinear parabolic
problems capable to be numerically solved in a highly efficient way. Finally,
in order to assess the computational feasibility of the algorithm, we have compared
our results with those obtained using competitive (freely available) parallel
numerical codes, which are widely used by the high-performance scientific
computing community.

Here it is the outline of the paper. In Sec. \ref{agpr} a generalized
probabilistic re\-presentation is presented, discussing two different possible 
strategies based on suitable random trees. Moreover, a qualitative study of 
the numerical errors is accomplished analyzing a few relevant test examples.
Sec. \ref{examples} is devoted to numerical examples, where the high efficiency of the PDD method 
comparing with classical methods is illustrated.
Finally, we summarize the more relevant findings to close the paper. 

\section{A generalized probabilistic representation}
\label{agpr}

In order to generalize the class of parabolic problems amenable to a probabilistic representation in terms of branching diffusion processes, it becomes more convenient to rewrite Eq.(\ref{nonlinear_new}) in an 
integral form. This can be done readily resorting to the Duhamel principle \cite{duchateau} for inhomogeneous initial-value parabolic problems, and
yields
\begin{equation}
 u(x,t) = \int_{{\bf R}^n} dy \, g(y) \, p(x,t,y,0) 
     + \int_0^t \int_{{\bf R}^n} ds \, dy \, f(u(y,s),y,s) \, p(x,t,y,s),\label{integraleq}  
\end{equation}
where $p(x,t,y,\tau)$ is the associated Green's function, satisfying the
equation
\begin{eqnarray}
    \frac{\partial p}{\partial t} = Lp ,     \   \  
           x \in {\bf R}^n,  \quad  t > \tau \nonumber  \\
       p(x,\tau,y,\tau) = \delta(x-y).
          \label{green2}
\end{eqnarray}

The main difference with the previous representation obtained in~
(\ref{nonlinear_prev}) 
rests on the absence of the constant potential term $-c u(x,t)$. Such a term 
was crucial, since it allowed to obtain a probabilistic representation based on generating {\it branching} 
diffusion processes governed by an exponential random time, S, with density probability $p(S)=c \,\exp(-c S)$ in \cite{Angel1,Angel2}. In the following we propose two different strategies capable to overcome such a constraint generalizing further the aforementioned representation. 

\subsection{Strategy A}

This first strategy consists in inserting artificially a constant potential term
into the PDE, by simply changing variables as follows
\begin{equation}
   u(x,t) = v(x,t) e^t. \nonumber
\end{equation}
Then, $v(x,t)$  should satisfy the following equation:
\begin{eqnarray}
   v(x,t) &=& e^{-t} \int_{{\bf R}^n} dy \, g(y) \, p(x,t,y,0) \nonumber    \\
     &+&  \sum_{j=2}^m\int_0^t \int_{{\bf R}^n} ds \, dy \, e^{- s} \, c_j(y,t-s) 
e^{(j-1)(t-s)}v^{j}(y,t-s)
 \, p(x,s,y,0).                                                        \label{totalduhamel}
\end{eqnarray}

The Green function can be obtained probabilistically 
 by means of the celebrated Feynman-Kac formula \cite{Karatzas} as follows
\begin{equation}
	p(x,t,y,\tau) = E [ \delta(\beta(t)-y)], \label{green}
\end{equation}
where $\beta(t)$ is the solution of an initial-value problem for the  
stochastic differential equation (SDE) of the Ito type, related to the elliptic
operator in (\ref{nonlinear_new}), i.e.,
\begin{equation}
        d \beta = b(\beta,t) \, dt + \sigma(\beta,t) \, d W(t), \quad  
\beta(\tau)=x. \label{SDE}
\end{equation}   
Here $W(t)$ represents the N-dimensional standard Brownian motion (also
called Wiener process); see \cite{Karatzas}, e.g., for generalities, 
and \cite{Kloeden,Milstein} for related numerical treatments. 
The drift, $b$, and the diffusion, $\sigma$, in (\ref{SDE}), are related to 
the coefficients of the elliptic operator in (\ref{nonlinear_new}) by $\sigma^2 
= a$, with $a > 0$.
Substituting (\ref{green}) into Eq.(\ref{totalduhamel}), and introducing 
 $\mathbf{1}_{[S > t]}$ as the indicator (or characteristic) 
function, being $1$ or $0$ depending whether $S$ is or is not 
greater than $t$,  Eq.(\ref{totalduhamel}) can be rewritten as follows

\begin{eqnarray}
\label{eq:probexp1}
 v(x,t) &=& E \left[ g(\beta(t)) \, \mathbf{1}_{[S>t]} \right] \\ \nonumber
 &+& E \left[ \sum_{j=2}^m c_j(\beta(S),t-S) \, e^{(t-S)(j-1)} \, v^{j}(\beta(S),t-S) \, \mathbf{1}_{[S \leq t]} \right].
\end{eqnarray}

Here the time $S$ is a random time, drawn from the exponential density
distribution $p(S)= \exp(-\, S)$. 

The equation above can be recursively solved, 
replacing the last term on 
the right-hand side with the solution $v(x,t)$, obtaining in such way 
an expansion in terms of multiple exponential random times, $S_i$, similarly as it was done in 
\cite{Angel1,Angel2}. However, rather here the procedure is much more involved
since now the integral equation contains both, variable coefficient terms,  $c_j(\beta(S),t-S) \, e^{(t-S)(j-1)}$, and various nonlinear terms labeled by $j$.   
The latter can be reformulated probabilistically introducing a new discrete random variable $\alpha$ taken values 
between $2$ and $m$, and governed by a 
uniform probability distribution with probability $p=1/(m-1)$, as follows

\begin{eqnarray}
\label{eq:probexp2}
 v(x,t) &=& E \left[ g(\beta(t)) \, \mathbf{1}_{[S>t]} \right] \\ \nonumber
 &+& E \left[ c'_\alpha(\beta(S),t-S) \, e^{(t-S)(\alpha-1)} \, v^{\alpha}
(\beta(S),t-S) \, \mathbf{1}_{[S \leq t]} \right],
\end{eqnarray}
where  $c'_i=(m-1)c_i$. Note that the probability distribution for the random variable $\alpha$ can be chosen in principle arbitrarily, however in order to minimize
the statistical error it turns out convenient to assume $\alpha$ uniformly distributed since in such way all the nonlinear terms are equally sampled. Expanding recursively the equation above, we obtain:
\begin{eqnarray}
&&v(x,t) = E \left[ g(\beta(t)) \, \mathbf{1}_{[S_0>t]} \right] \nonumber\\
 && + E\left[h^{(\alpha_1)}(y_1(S_0),S_0)\prod_{i=1}^{\alpha_1} g(x_i(t-S_0)) \, \mathbf{1}_{[S_i>t-S_0]} \, \mathbf{1}_{[S_0 < t]}\right] \nonumber\\
&& +E\left[h^{(\alpha_1)}(y_1(S_0),S_0)
h^{(\alpha_2)}(y_2(S_0+S_1),S_0+S_1)  
\mathbf{1}_{[S_i>t-S_1]}
\right. \nonumber\\
&& \times \left. \prod_{i=2}^{\alpha_1} g(x_i(t-S_0)) \prod_{i=\alpha_1+1}^{\alpha_1+\alpha_2} g(x_i(t-S_0-S_1)) 
\mathbf{1}_{[S_i>t-S_0-S_1]}
\mathbf{1}_{[S_1<t-S_0]}
\mathbf{1}_{[S_0<t]}\right] \nonumber\\
&&+\cdots,\label{recursion}
\end{eqnarray}
where $h^{(\alpha)}(y,s)=c'_\alpha(y,t-s) \, e^{(t-S)(\alpha-1)}$. Here 
$x_i(t)$, $y_i(t)$ denote the position of the $i$th path of the stochastic 
process $\beta(t)$ surviving or expiring, respectively,  at time $t$.

Note that in Eq.~(\ref{recursion}) the solution can be evaluated by simply summing each partial contribution, and similarly to the class of equations studied in \cite{Angel1,Angel2}, the computational tool based on generating random trees turns out to be very useful since it allows 
to rapidly obtain analytically or compute numerically such contributions.

In the following, we assume the reader is familiar with 
the terminology used here, like the concepts of \textit{branch}, \textit{random tree}, etc.
However, it is useful at this point to recall some terminology pertaining to 
trees, usually directed trees. A tree is a connected graph, i.e., a set of 
nodes linked by {\it edges}, with only one starting node, called {\it root}, a 
number of final nodes, called {\it leaves}, while the other (internal) nodes 
are called {\it vertices}, and we call here {\it branch} every set of edges 
joining vertices to leaves. Therefore, in this sense 
the number of leaves in a tree is equivalent to the number of branches. Finally, the number of nodes linked to a given node is called the number of {\it children}. 

\begin{figure}[ht] 
\begin{center}
\includegraphics[width=0.9\textwidth,angle=0]{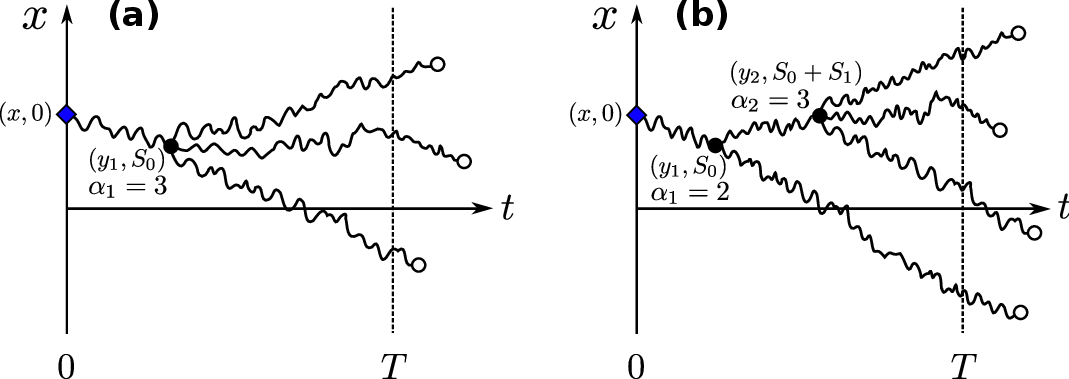}
\caption{Two possible configuration diagrams. In (a) only one splitting event occurs, with $\alpha_1=3$. Rather in (b) two of such events occur, the first one with $\alpha_1=2$, and the second one with $\alpha_2=3$.} 
\label{fig:bbm_two_scenarios}
\end{center} 
\end{figure}

 The algorithm works as follows: We first generate a random 
exponentially distributed time, $S_0$, and a random path belonging
to the stochastic process $\beta$. If $S_0$ is less
than the final time $t$, then we split the given path into as many branches
as those corresponding to the randomly chosen value $\alpha$, the degree of nonlinearity. They depart 
from the position where the previous path was at time $S_0$, and continue
along independent trajectories until the next splitting event takes 
place. Whenever one of the possible branches reaches the final time, $t$, the 
initial value, $g$, is evaluated at the position where the path 
was located. Finally, the partial contribution to the solution is reconstructed multiplying all 
contributions coming from each branch and the coefficients $h^{(\alpha)}(y,s)$
conveniently evaluated at specific points. For the purpose of illustration, we sketch a typical configuration in the second picture of 
Fig.~\ref{fig:bbm_two_scenarios}, corresponding to two different splitting events. Note that 
such a configuration represents graphically what appears in the last term of 
Eq.~(\ref{recursion}). Moreover, it can be seen clearly the structure of the random tree, where the dots denotes the splitting events, being marked in black or white depending on whether the splitting time obtained is less or not than the final time, $T$. Therefore, the white dots correspond to the leaves of the tree, while
the black dots are the corresponding vertices. Note that the number
of children associated to a given node depends on
the value of $\alpha$, which has been randomly chosen.

Hence the solution $v(x,t)$ in Eq. (\ref{eq:probexp2}) can 
be rewritten as a  
expectation value over random trees of a suitable multiplicative functional of the initial data $g(x)$, the random times $S_i$, and the variable coefficients  $c'_j$, as follows,
\begin{eqnarray}
v(x,t)&=& E \left[\prod_{i=1}^{Ne(\omega)}h^{(\alpha_i(\omega))}(y_{i}(\omega),
\bar{S}_{i}(\omega))  
\prod_{l=1}^{k(\omega)}g(x_l)  \right] \label{rep_strategyB}.
\end{eqnarray}
Here $k(\omega)$, and $Ne(\omega)$ are the random number of branches at final time $t$, and the number of splitting events obtained when generating the random tree, respectively. By $\bar{S}$ we denote the corresponding global random time obtained by
summing conveniently the random times $S_i$ according to the specific structure of the generated random tree. It is worth to observe that such trees are used as a tool to construct the structure representing a given partial 
contribution to the solution, allowing afterward to follow easily how the arguments of the
functions $h^{(\alpha)}(y,s)$ are
exchanged when solving recursively Eq.(\ref{eq:probexp2}). 

Even though such a new representation allows to expand further the class of 
equations suited to be computed probabilistically, however
rather than in the representation obtained  in \cite{Angel1,Angel2}, a major drawback now should be faced. This is because the coefficients multiplying the nonlinear terms $v^j$ might be greater than $1$, therefore being the convergence of the numerical procedure not guaranteed. In fact, the series obtained by 
expanding Eq. (\ref{eq:probexp2}) could be divergent, and in general cannot be summed simply by a sequence of partial sums. Additionally, the \textit{pruning} techniques presented in \cite{Angel1,Angel2} cannot be applied for the same reason. Nevertheless, numerical experiments
show that in many cases the asymptotic series can be summed up resorting to summation methods such as the Euler's formula, or approximation techniques based on the Pade approximant. In this paper, we consider
merely the latter one, because from numerical experiments seems to be more robust in dealing with the unavoidable numerical errors affecting the coefficients of the series computed numerically. In Sec. \ref{qsotne} a few test problems have been
investigated to illustrate the robustness and convergence of the Pade approximant.

\subsection{Strategy B}

A different strategy for obtaining a probabilistic representation for the problem in Eq.~(\ref{nonlinear_new}) consists in sampling both terms of the integral equation (\ref{integraleq}), by introducing a two-point discrete random variable $\xi$ taking the values 
$0$, and $1$, with probability 
$P(0)=q,\, P(1)=1-q$. 
Therefore, the integral equation (\ref{integraleq}) can be rewritten as follows,

\begin{eqnarray}
\label{totalduhamel_B}
   u(x,t) &=& q \, \int_{\Omega} dy \, \tilde{g}(y) \, p(x,t,y,0) \nonumber    \\
     &+& (1-q)\int_0^t \int_{\Omega} ds \, dy \, \sum_{j=2}^m \tilde{c}_j(y,t-s) u^{j}(y,t-s) \, p(x,s,y,0),
\end{eqnarray}
where $ \tilde{g}(x)=g(x)/q$, and $\tilde{c}_j(x,t)=c_j(x,t)/(1-q)$. The
probabilistic representation can be readily found and has the form 
\begin{eqnarray}
\label{eq:probexpB1}
 u(x,t) &=& E \left[\tilde{g} (\beta(t))\delta(\xi)\right] \nonumber\\
&+& E \left[\eta(t) \tilde{c}'_{\alpha}(\beta(tS),t(1-S)) \, u^{\alpha}(\beta(tS),t(1-S))\delta(\xi-1)\right],
\end{eqnarray}
where the time $S$ is a random time between $0$ and $1$ uniformly distributed,
$\alpha$ a discrete random variable taking values 
between $2$ and $m$ with equal probability $p=1/(m-1)$, $\tilde{c}'_{\alpha}=(m-1)\tilde{c}_{\alpha}$, and $\eta(t)=t$. Similarly to the previous strategy, the equation above can be recursively expanded, and yields

\begin{eqnarray}
\label{eq:probexpB2}
&&u(x,t) = E \left[ \tilde{g}(\beta(t)) \, \delta(\xi_0) \right]\nonumber\\
&&+E\left[\eta(t)\tilde{c}'_{\alpha_1}(y_1(tS_0),t (1-S_0))\prod_{i=1}^{\alpha_1} 
\tilde{g}(x_i(t(1-S_0)) \, \delta(\xi_i)\delta(\xi_0-1)\right]
\nonumber\\
&&+E\left[\eta(t)\eta(t (1-S_0))\tilde{c}'_{\alpha_1}(y_1(t S_0),t(1-S_0))
\right. \nonumber\\
&&\times \, \tilde{c}'_{\alpha_2}(y_2(t (1-S_0)S_1),t(1-S_0)(1-S_1)) \prod_{i=2}^{\alpha_1} g(x_i(t(1-S_0)))) \, 
\delta(\xi_i) \nonumber\\
&&\left. \times \prod_{j=\alpha_1+1}^{\alpha_1+\alpha_2+1} g(x_j(t(1-S_0)(1-S_1))) 
\delta(\xi_j)\delta(\xi_1-1)\delta(\xi_0-1)\right] +\cdots
\end{eqnarray}

Therefore, as in the strategy A, the solution can be obtained as the expectation
value over suitable random trees of a given multiplicative functional of 
the initial condition, being given now as follows,
\begin{eqnarray}
\label{eq:probexpBclosed}
u(x,t)&=& E \left[\prod_{i=1}^{Ne(\omega)}\eta(t\bar{S}_{i}(\omega))c_{\alpha_i(\omega)}(y_{i}(\omega),
t\bar{S}_{i}(\omega))  
\prod_{l=1}^{k(\omega)}g(x_l)  \right].
\end{eqnarray}
Here $\bar{S}$ is the global time random variable associated to the generated random tree. It is obtained multiplying the different random times $S_i$ according to the specific structure of the tree. 

In view of the probabilistic representation obtained for both strategies, 
it is worth to point out that both strategies require generating random trees to evaluate numerically the partial contribution to the solution, however in practice the computational procedure needed is quite different. In fact, while the random trees in the strategy A are constructed by generating a unique random number, the random time $S$,  which governs how the trees branch off in time, rather those in the strategy B require two independent random numbers for the same purpose. This is because in the strategy A the change of variable introduces a time dependent exponential coefficient, which can be used to construct a probabilistic representation for Eq.~(\ref{eq:probexp1}) based on an exponential random time, and therefore the random trees can be fully characterized by such random time. On the contrary, the probabilistic representation in Eq.~(\ref{eq:probexpB1}) requires both, a random number $\xi$ which governs the branching process, and a random time $S$ uniformly distributed, for evaluating numerically the partial contribution to the solution corresponding to a given random tree.

Similarly to the strategy A, the series obtained using the strategy B in Eq. 
(\ref{eq:probexpB2}) may be divergent, being therefore necessary to resort to approximation techniques, such as the Pade approximant, to approximate conveniently the sum of the series.

To illustrate how both strategies can be implemented in practice for solving
an initial value semilinear parabolic problem, let consider the following equation,
\begin{eqnarray}
  && \frac{\partial u}{\partial t} = \frac{\partial^2 u}{\partial x^2} + u^2, 
      \quad  x \in {\bf R},  t > 0  \\
   && u(x,0) = f(x).              \label{problem1}
\end{eqnarray}
Since the procedure underlying the strategy A is formally similar
to that followed when a constant potential term is present, we refer to the reader to \cite{Angel1,Angel2} for a practical implementation for such a case, and here we focus merely on the strategy B. From Eq. (\ref{eq:probexpB1}), the probabilistic representation is given by
\begin{eqnarray}
 u(x,t) = E \left[\tilde{g} (\beta(t))\delta(\xi)\right] 
+ E \left[\eta(t) \, u^{2}(\beta(tS),t(1-S))\delta(\xi-1)\right],
\label{probrep_coin_problem1}
\end{eqnarray}
or in a more compact format, using Eq.(\ref{eq:probexpBclosed}) for the expectation value over random trees of a given multiplicative functional in Eq.(\ref{probrep_coin_problem1}), by
\begin{eqnarray}
u(x,t)&=& E \left[\prod_{i=1}^{Ne(\omega)}\eta(t\bar{S}_{i}(\omega))
\prod_{l=1}^{k(\omega)}g(x_l)  \right].
\end{eqnarray}
Every random tree is built generating a sequence of interconnected binary random variables, $\xi_i$, branching off from the previous one as follows: Let $\xi_1$ the random variable associated to the root of the tree. Only when $\xi_1$ takes value 1 with probability $P(1)=1-q$, two new random variables denoted by $\xi_{2,3}$ (child nodes of the root), are created. These new variables proceed further creating other nodes governed by the same rule, until no random number $\xi_i$ takes anymore the value 1. At this point the procedure is concluded, giving rise to a random tree characterized by $k$ branches or leaves, and $Ne$ splitting events. 

The nodes of the tree are labeled in binary format according to their ancestors as follows: A given node with label $[a_0a_1a_2...a_N]$, where $a_i=0,1$, is connected to the set of nodes $\{[a_0],[a_0a_1], [a_0a_1a_2],\ldots,[a_0a_1a_2\cdots a_{N-1}]\}$. The global time random variable $\bar{S}$ associated to a given tree with $k$ branches is given by
\begin{eqnarray}
\bar{S}=\prod_{i=1}^{2^{k-1}-1}S_j^{\gamma_j}, \quad 
\gamma_l=\sum_{j=l+1}^{2^{k-1}-1}\nu_j\left<j|l\right>,\,l=1,\ldots,2^{k-1}-1
\end{eqnarray}
where $\nu_l$ is 0, or $1$ depending on whether the tree contains or not the node $l$. The function $\left<\cdot|\cdot\right>$ is defined as follows,
\begin{eqnarray}
\left<j|l\right>:= \left\{ \begin{array}{ll}
         1 & \mbox{if $T_j^{[l]}=l$}\\
        0 & \mbox{otherwise}.\end{array} \right.,
\end{eqnarray}

where both, $j$ and $l$ are numbers written in binary format, and $T_j^{[l]}$ is an operator that truncates the number $j$ to their most significant $[l]$ digits,
where $[l]$ is the number of digits of $l$. By example, let $j=[a_0a_1a_2...a_N]$, then  $T_j^{[l]}=[a_0a_1...a_{[l]-1}]$.  

\begin{figure}[ht] 
\begin{center}
\includegraphics[width=0.8\textwidth,angle=0]{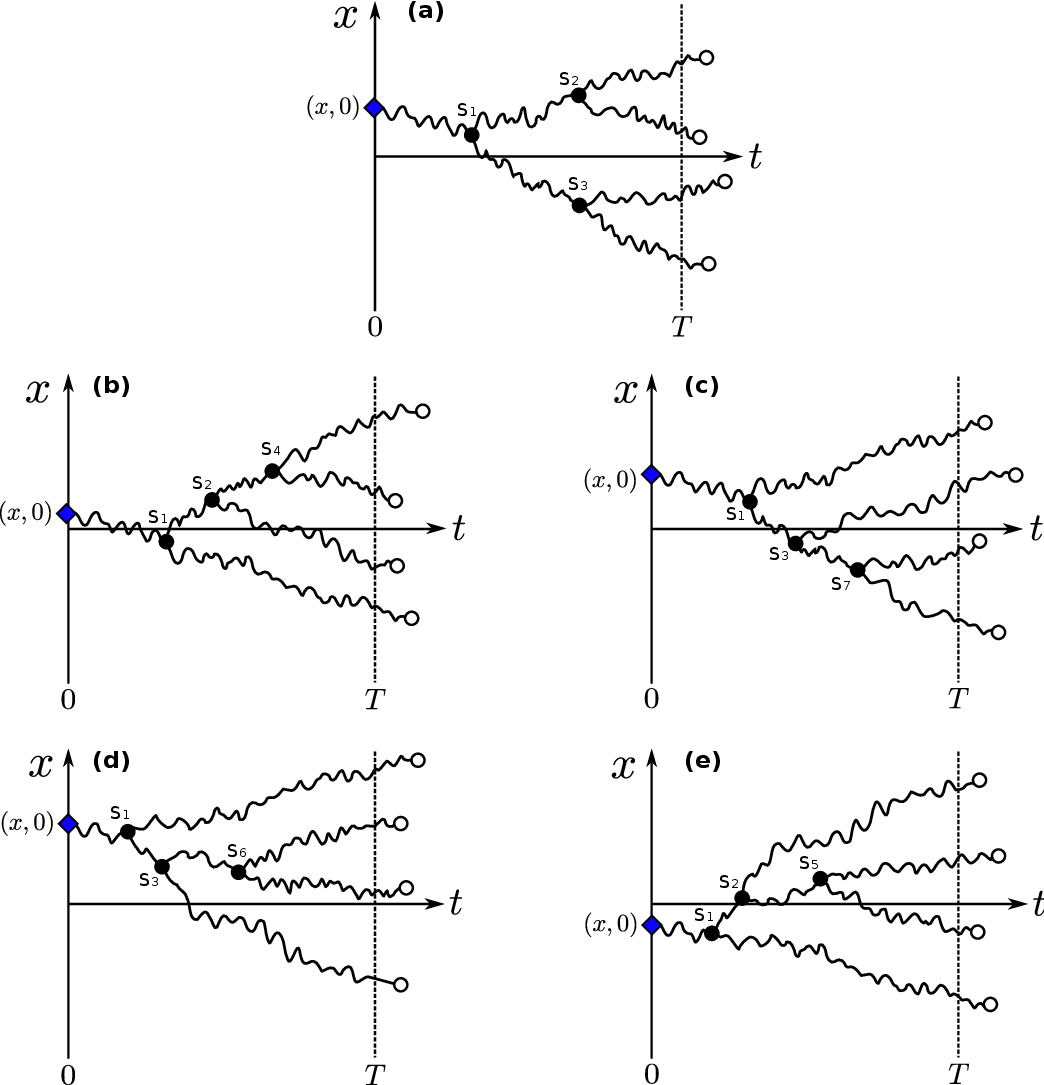}
\caption{Configuration diagram for the case of 4 branches and 3 splitting events. Here $S_i$ is a random time uniformly distributed between the previous generated time, and the final time, $T$. The 
corresponding labels $i$ of the random time $S_i$ are defined according to the rule explained in the text.} 
\label{fig:bbm_k4_Ne3}
\end{center} 
\end{figure} 

Figure \ref{fig:bbm_k4_Ne3} shows the different random trees obtained with $k=4$, and $Ne=3$, and their corresponding labels according to the rule 
defined above.


\subsection{Computational complexity of the strategy B}

  In this subsection we estimate the computational complexity in terms of
the computational time required to compute probabilistically the 
solution at a 
single point $(x,t)$  based solely on the strategy B, since the computational complexity of the strategy A coincides with that analyzed already in \cite{Angel1,Angel2}.

  The branching stochastic process associated to the nonlinear term $u^{j}$, 
  requires creating $j$ branches every time a splitting event occurs,
being therefore more costly whenever the power of the nonlinearity is higher. 
Then, it is worth to observe that for the general nonlinear function in Eq.~(\ref{nonlinear_new}), the overall computational time is governed by the computational time spent by the nonlinear coefficient with the highest power, say $u^m$.


  The computational time spent to generate any given branch, is a function of 
the final time, $t$, as well as of the time step, $\Delta t$, chosen to 
solve numerically the associated stochastic differential equation in 
(\ref{SDE}). In addition, we should take into account the random times 
$\Delta t_s$ responsible for branching. From Eq. (\ref{eq:probexpB1}), it holds that
$\Delta t_s=t\,S$, where $S$ is a random number picked up from 
the uniform distribution $U(0,1)$. It is thus necessary 
that the time-step discretization, used to solve (\ref{SDE}), also captures 
the instants when the random exponential times occur. In practice, the actual
time step is chosen according to the minimum value between $\Delta t$ and 
$\Delta t_s$. 
Since $\Delta t_s$ is chosen randomly, the probability of being less 
than $\Delta t$ can be easily estimated, and turns out to be $(\Delta t)/ t$. When this occurs, the actual time 
step used for the numerical solution of (\ref{SDE}) should be chosen to be
$\Delta t_s$. 
Averaging over all random trees, we obtain the most probable 
time step to be used, which is given by
\begin{eqnarray} 
   \overline{\Delta t_s} = 
            \int_0^{\Delta t} ds \,s \, \frac{1}{t}+\Delta t 
 \int_{\Delta t}^t ds \,\frac{1}{t}
         = \Delta t - \frac{1}{2} \frac{(\Delta t)^2}{t}.
\end{eqnarray}   

  The computational time can be measured, typically, in terms of the number of 
iterations in time, required to fully generate a random tree 
with $k$ branches up to the final time, $t$. Defining $t_c$ as the time spent 
per iteration, such computational time can be estimated as $k t_c \, 
t/\overline{\Delta t_s}$. In case of $N$ random trees, the average 
computational time, $t_b$ ($b$ standing for ``branching''), turns out to be 
\begin{equation}
  t_b = N \sum_{k=1}^\infty k t_c {\frac{t}{\overline{\Delta t_s}} P(k,m)}, 
                                                                 \label{Tb}
\end{equation}
where $P(k,m)$ is the probability of finding a random tree 
with $k$ branches, being $m$ the number of children.
 
\begin{figure}[ht] 
\begin{center}
\includegraphics[width=0.8\textwidth,angle=0]{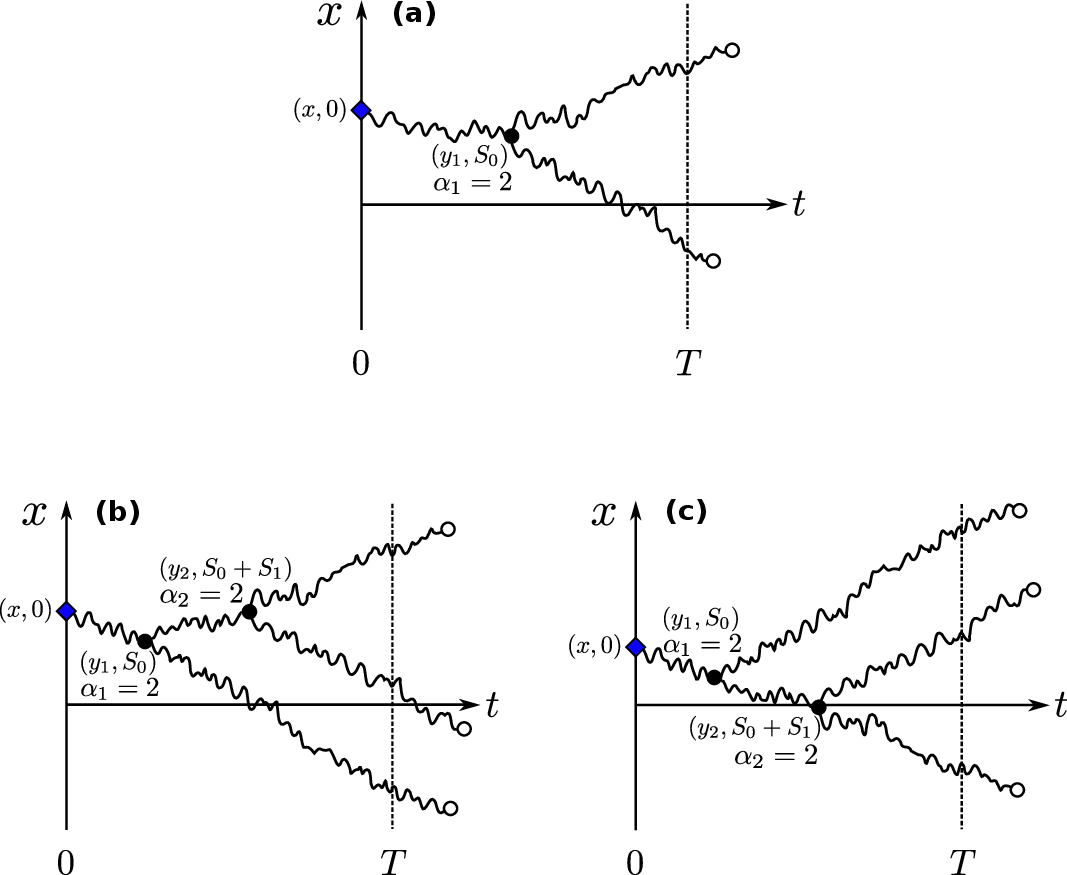}
\caption{Configuration diagrams for the case of 2, and 3 branches, illustrating the notation used in the probabilistic representation obtained for Strategy A in Eq. (\ref{rep_strategyB}). 
Here $\alpha_i$ denotes the 
power of the nonlinearity, chosen randomly between $2$ and $m$, and governing the number of branches 
created every time a splitting event $i$ takes place; $y_i$ denotes the position of the $i$th path of the stochastic process expiring at a given time ${\bar S}$. By ${\bar S}$ we denote the corresponding global time obtained by summing conveniently the random times $S_i$ according to the rules described in the text. } 
\label{fig:bbm_k234}
\end{center} 
\end{figure}

\begin{figure}[ht] 
\begin{center}
\includegraphics[width=0.5\textwidth,angle=0]{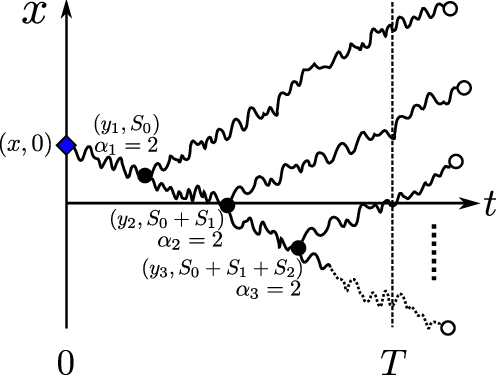}
\caption{Configuration diagram with $k$ branches. See Fig. 3 for a detailed explanation of the 
corresponding labels.} 
\label{kbranches}
\end{center} 
\end{figure}

  Such a probability can be evaluated by first enumerating and then summing up 
the various probabilities,  $p_i(k)$, of having $k$ branches as a final 
configuration, that is
\begin{equation}
          P(k,m) = \sum_{i=1}^{N_D(k,m)} p_i(k,m).                       \label{pk}
\end{equation}
Here $N_D(k,m)$ denotes the total number of possible diagrams characterized 
by $k$ branches and $m$ children, and $p_i(k,m)$ the probability of obtaining each of them.
In Fig.~\ref{fig:bbm_k234} we show, for the purpose of illustration, some 
diagrams for $k = 2, 3$. In particular, for $k=4$ the total number 
of possibilities of obtaining $4$ branches has been shown in Fig. 
\ref{fig:bbm_k4_Ne3}. It is reasonable to assume 
that each diagram contributes equally to the probability function (\ref{pk}). 
Therefore, such a probability can be obtained by simply counting the number of 
possible diagrams with $k$ branches, $N_D(k,m)$, and then multiplying by the 
probability of having one of them, that is $P(k,m) = N_D(k,m) p_1(k,m)$. For 
convenience, we consider the special diagram shown in Fig.~\ref{kbranches}. 
The probability of obtaining such a diagram as a final configuration is given 
by

\begin{eqnarray}
   p_1(k) = q^k (1-q)^{Ne}.        \quad \quad
\end{eqnarray}
Note that the number of branches, $k$, is related to 
the number of splitting events, $N_e$, by $k = (m-1)N_e + 1$. 

Concerning the number of possible diagrams, $N_D$, this should be a function
of the number of branches, $k$.  A simple strategy to compute the all possible configurations requires analyzing the distribution of leaves at any depth level of the random tree.

For simplicity, in the following let consider first the case of $m=2$, 
which corresponds to binary trees, and a tree composed of $k$ branches.
When the tree is generated, the first splitting event gives rise to two new 
sub-trees, each one having at least one leaf, and at most $k-1$ leaves, since the total number of leaves should be $k$. The different possible ways of distributing the $k$ leaves among the two sub-trees determines the different possible configurations at this level, namely let assign $1$ leaf to the first sub-tree, and $k-1$ to the second one; or $2$ to the first one, and $k-2$ to the second one, and so on. Clearly, this procedure should be recursively applied to every sub-tree, since any of them exhibits different type of configurations. As a result the following nonlinear, full memory, convoluted recurrence is obtained
\begin{eqnarray*}
       N_{D}(k,2) &=& \sum_{j=1}^{k-1}  N_{D}(k-j,2) \ N_{D}(j,2),\\
       N_{D}(0,2) &=& 0,\quad N_{D}(1,2) = 1,  \label{recursive}
\end{eqnarray*}
see~\cite{Graham} e.g. This recurrence can be solved by means of generating functions, yielding
 \begin{eqnarray*}
     N_{D}(k,2) &=&\frac{1}{k} \left( {\begin{array}{*{20}c} 2k - 2 \\ k-1 \\ \end{array}} \right).
                                                           \label{multif}
\end{eqnarray*}

Moreover, this can be generalized further to any number of children by considering the Fuss-Catalan numbers, see~\cite{Aval}. Recall that the number of 
children is given by the power of the nonlinearity $m$ in Eq.(\ref{nonlinear_new}), and therefore the corresponding tree should be in general an $m$-ary tree.  For this case it becomes more convenient to describe the number of configurations $N_D$ in terms of the number of splitting events $Ne$ instead of using the number of branches $k$, since when $m$ is different from 2, $k$ is not a consecutive integer number. The solution is given by

\begin{eqnarray*}
     N_{D}(Ne,m) = \frac{1}{Ne \ m + 1} \left( {\begin{array}{*{20}c} Ne \ m + 1 \\ Ne \\ \end{array}} \right).
                                                           \label{multif}
\end{eqnarray*}

  Summarizing, the probability of obtaining a random tree with $k$ branches 
is given by 
\begin{eqnarray}
    P(k,m) = q^k (1-q)^{Ne}  \frac{1}{m\,Ne+1} \left( {\begin{array}{*{20}c} 
m\,Ne +1\\ Ne \\ \end{array}} \right).
                                                          \label{pfinal}
\end{eqnarray}

However, it turns out that the function $P(k,m)$ in Eq.(\ref{pfinal}) can be considered a probability function for a particular range of values of $q$, since only for such values $\sum_{Ne=0}^{\infty} P(k,m)=1$ is satisfied, being different from one or even unbounded otherwise. In fact, let $u=\sum_{Ne=0}^{\infty} P(k,m)$, $v=u/q$, and $x=q^{m-1}(1-q)$, then

\begin{equation}
\label{eq:binomial_series}
v=\sum_{Ne=0}^{\infty} x^{Ne}  \frac{1}{m\,Ne+1} \left( {\begin{array}{*{20}c} 
m\,Ne +1\\ Ne \\ \end{array}} \right).
\end{equation}
Note that the power series above coincides with the generalized binomial series \cite{Graham}, thus $v=B_m(x)$. From the properties of the generalized binomial series, it holds that
\begin{equation}
\label{eq:binomial_series2}
v=1+x\,v^m.
\end{equation}
Among the $m$ possible solutions of Eq.(\ref{eq:binomial_series2}), only that one satisfying $\lim_{x\to 0} v=1$ is meaningful, which corresponds to the trivial probability function $P(k)=\delta_{k0}$ being $q=1$. The solution $v$ can be constructed iteratively applying a Picard iteration as follows,
\begin{equation}
v_n=1+x\,v_{n-1}^m,
\end{equation}
with $v_0$ arbitrarily chosen. Note that this can be seen as a dynamical map
$v_{n}=f(v_{n-1})$, with $f(v)=1+x\,v^m$. Such a map has as a fixed point $v^{*}=1/q$, being stable whenever $f'(v^{*})\leq 1$. This corresponds to values of $q$ satisfying 
\begin{equation}
q\geq (m-1)/m.\label{q}
\end{equation}
This can be generalized further for a nonlinear function as in Eq.~(\ref{nonlinear_new}),
obtaining the range of allowed values of $q$ for which a probability function
$P(k,m)$ can be found. Let define now the function 
$u=\sum_{k=1}^{\infty}\sum_{l=\lceil (k-1)/(m-1)\rceil}^{k-1} P(k,l)/(m-1)$
Similarly to the previous case, it can be readily
proved that $u$ satisfies the following equation,
\begin{eqnarray}
u=q+\frac{(1-q)}{m-1}\sum_{i=2}^{m} u^{i}.
\end{eqnarray}
The solution $u=1$ can be iteratively constructed by applying a Picard iteration, thus obtaining the following nonlinear map,
\begin{eqnarray}
u_n=f(u_{n-1}),
\end{eqnarray}
where $f(u)=q+\frac{(1-q)}{m-1}\sum_{i=2}^{m} u^{i}$.
Since $f'(1)=(1-q)[1+(m+1)/2]$, the fixed point $u=1$ turns out to be stable 
when 
\begin{equation}
q\geq m/(m+2).\label{qgen}
\end{equation}

\begin{figure}
  \centering
  \subfloat{\label{fig:pk_SB_m2}\includegraphics[width=0.34\textwidth,angle=-90]{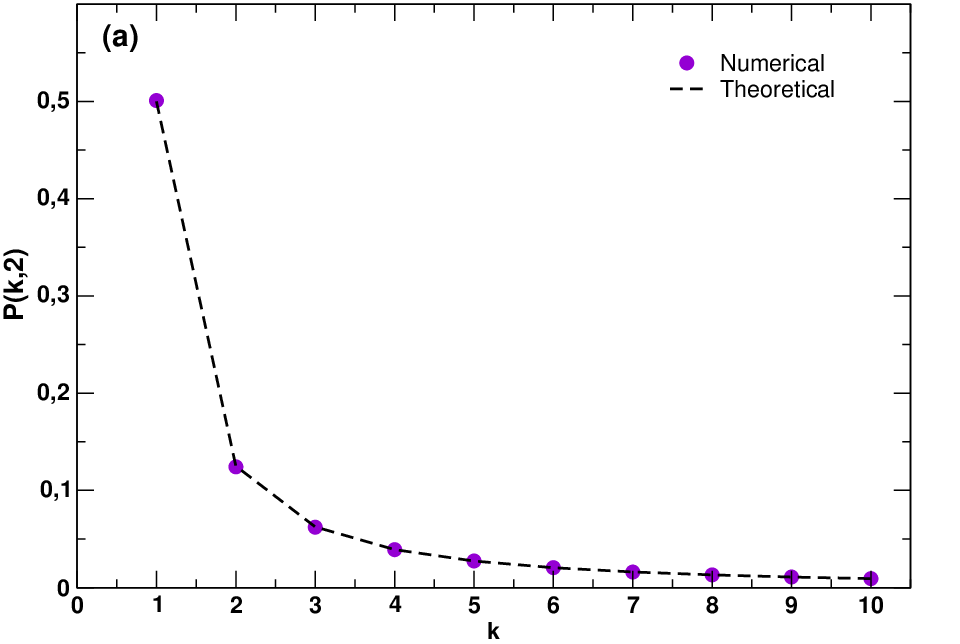}}
  \subfloat{\label{fig:pk_SB_m3}\includegraphics[width=0.34\textwidth,angle=-90]{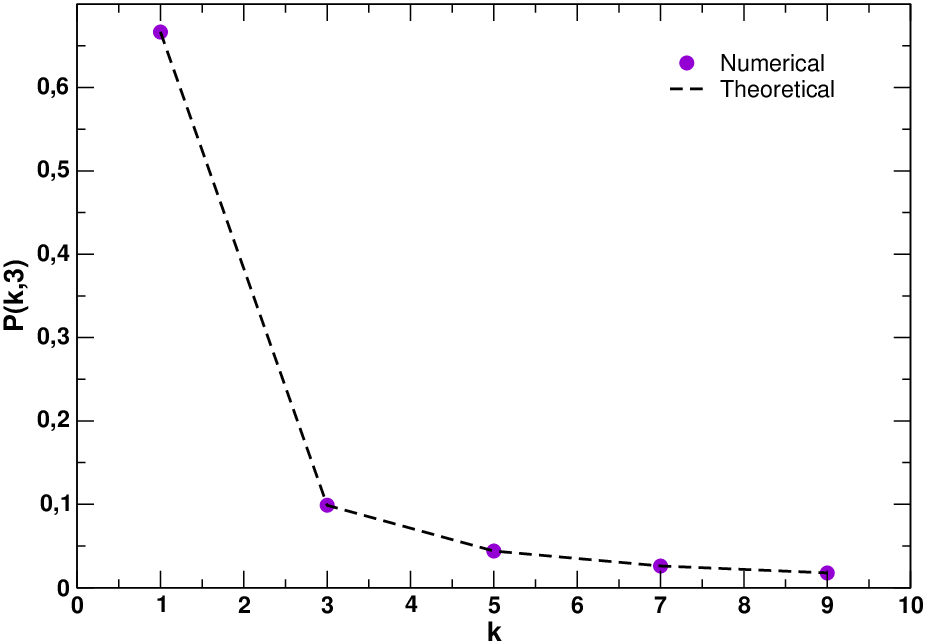}}
  \caption{Comparison between the probability function $P(k,m)$ obtained
analytically and numerically simulating $10^6$  random trees for $m=2$ in (a), and $m=3$ in (b).}
  \label{pkcompare}
\end{figure}

The validity of Eq. (\ref{pfinal}) can be confirmed by some numerical
simulations, consisting in generating $N$ random trees with a given $q$ satisfying the specific constraint described above, and counting the number of branches obtained. A comparison between the probability $P(k,m)$ as function of $k$, obtained both 
numerically and theoretically, is plotted in Fig.~\ref{pkcompare}. This has 
been done for the case $m=2$ in Fig.~\ref{pkcompare}(a), and $m=3$ in 
Fig.~\ref{pkcompare}(b).
The perfect agreement validates the formula (\ref{pfinal}).

  Once the probability function, $P(k,m)$ is known, the computational time 
$t_b$ spent by the probabilistic part of the algorithm, can be evaluated. 
 From (\ref{Tb}), we have
\begin{equation}
    t_b \leq N t_c \frac{t}{\overline{\Delta t_s}} 
            \left<k\right>_{P(k,m)},                    \label{aproxTb}
\end{equation}
where $\left<k\right>_{P(k,m)}$ denotes the mean number of leaves, that is $\sum_{k=1}^{\infty} k P(k,m)$. Such a number can be easily computed exploiting the following
relation, obtained from Eq.(\ref{eq:binomial_series}) simply deriving with respect to $x$ and then
multiplying by $x$,
\begin{equation}
\sum_{Ne=0}^{\infty} {Ne}\, x^{Ne}  \frac{1}{m\,Ne+1} \left( {\begin{array}{*{20}c} 
m\,Ne +1\\ Ne \\ \end{array}} \right)=x\,\frac{dv}{dx}.
\end{equation}
Since for the allowed values of $q$, $v$ is given by $1/q$, then from 
Eq. (\ref{eq:binomial_series2}) $\left<k\right>_{P(k,m)}$ can be readily obtained and is given by
\begin{equation}
\left<k\right>_{P(k,m)}=\frac{q}{1-m(1-q)}.
\end{equation}
Hence an estimate of the computational time $t_b$ satisfies the following 
bound,
\begin{equation}
    t_b \leq N t_c \frac{t}{\overline{\Delta t_s}} \frac{q}{1-m(1-q)}
            .                    \label{aproxTb2}
\end{equation}
Note that the estimate obtained above exhibits a linear growth on $t$. This contrasts with the theoretical estimates of the computational time obtained for the strategy A \cite{Angel1,Angel2}, which grows instead unboundedly in time. Such a remarkable feature of this strategy in comparison with the strategy A allows to speed up notably the simulations. Moreover,
the bound obtained for the computational time suggests that decreasing $q$ may decrease further such a time, becoming singular however
when $q=(m-1)/m$, and therefore useless when such an occurrence happens.

\begin{figure}[ht] 
\begin{center}
\includegraphics[width=0.8\textwidth,angle=-90]{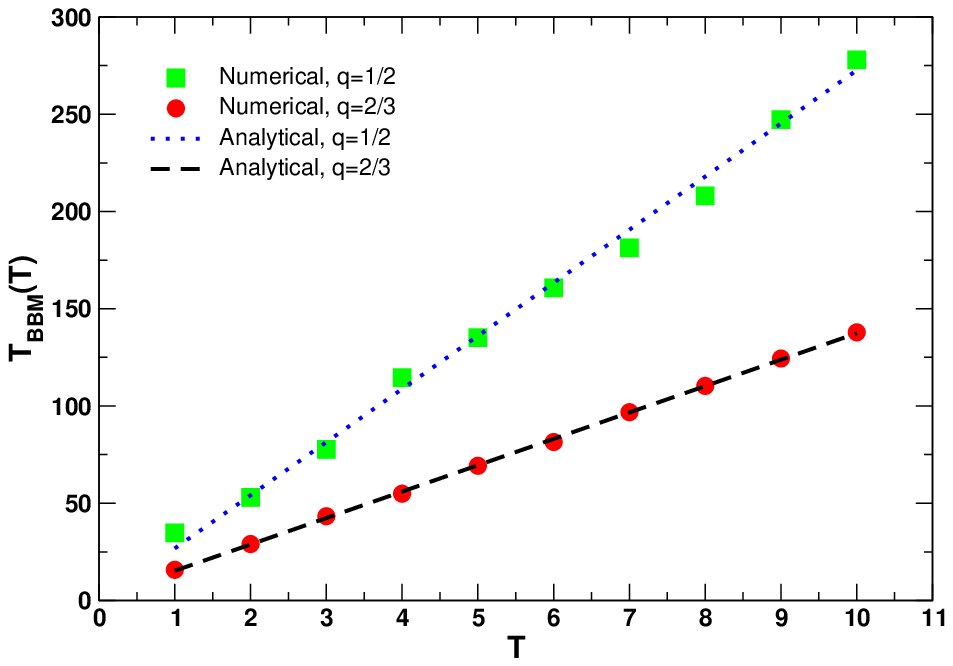}
\caption{Comparison between the computational time spent in solving 
Eq.~(\ref{problem1})  at a single point $x=0$, and the estimated computational time obtained theoretically in (\ref{aproxTb2}) for $m=2$. This has been done for $q=1/2$, and $q=2/3$. The numerical results have been fitted to a line with correlation coefficients $r=0.9973364$, and  $r=0.9998173$, respectively. Parameters are $N=10^6$.} 
\label{eq:tb_2_3} 
\end{center}
\end{figure}

  In Fig.~\ref{eq:tb_2_3}, a comparison between the theoretical estimates
obtained for $m=2$, and two different values of $q$,  and the measured computational times are shown as function
of the final time. Note the good agreement with the theoretical results.

%


\subsection{Qualitative study of the numerical errors}
\label{qsotne}

The numerical errors appearing when solving parabolic problems by some of the probabilistic strategies described above, are essentially the same as those analyzed in \cite{Angel1,Angel2}. In fact, the most
important source of numerical errors arises from replacing the expected 
value in Eq.(\ref{eq:probexpBclosed}) by a finite 
sum, and when the stochastic paths are actually simulated resorting to 
suitable numerical schemes. Thus, approximately, 
\begin{equation}
     u({\bf x},t) \approx \frac{1}{N}\sum_{j=1}^N f(\beta^*_j(t)),
\label{simplecase_discrete}
\end{equation} 
where $f$ denotes the multiplicative functional in Eq.(\ref{eq:probexpBclosed}), $N$ is the sample size, and $\beta^*$ is the stochastic path with 
discretized time. Clearly, such a discretization procedure unavoidably 
introduces two sources of numerical error. The first one  is the pure 
Monte Carlo statistical error, which it is known to be of order $O(1/\sqrt{N})$ when $N$ goes to infinity. The second error is due to the fact that the ideal 
stochastic path, $\beta_j(\cdot)$, has to be approximated, discretizing time, 
by some numerical scheme yielding the paths $\beta^*_j(\cdot)$. The truncation 
error made when solving numerically the stochastic differential equation 
(\ref{SDE}), obviously depends on the specific scheme chosen, see 
\cite{Kloeden}, e.g. Among these are the Euler scheme, which was used here to 
simulate numerically Eq. (\ref{SDE}). Such scheme is well known to have
a truncation error of order $O(\Delta t^{\alpha})$, where $\alpha=1/2$ or 
$\alpha=1$ depending on whether the scheme being of the ``strong'' or ``weak'' 
type, respectively \cite{Kloeden}.

Concerning the first error for the strategy B, the freedom to choose the value
of $q$ from a set of allowed values can be exploited in order to minimize such
an error. In fact, rather than the strategy A, the strategy B requires choosing specifically a given value of $q$ satisfying the constraints mentioned in the previous section. Among the set of allowed values, in the following we show that there exists an optimal one such that the statistical error made in computing probabilistically the coefficients of the Pade approximant turns out to be minimum. 

Clearly, the statistical error becomes larger when computing the contributions coming from trees with arbitrarily large number of branches, since for such a case
the number of generated trees is expected to be smaller, and consequently the associated 
multidimensional integral in (\ref{totalduhamel_B}) to be computed may be affected by a large statistical error. The number of generated trees $N_k$ with $k$ branches can be readily obtained, known the probability distribution $P(k,m)$, and is given by $N_k=N P(k,m)$, where $N$ is the sample size. Therefore the statistical error can be estimated, and it turns out to be of order of $O(N_k^{-1/2})$.

For simplicity let consider first the case of a nonlinear function in (\ref{nonlinear_new}) with a single nonlinear term, $u^j$. The value of $q$ minimizing the statistical error is attained when the number of random trees $N_k$ is maximum. Such a value can be obtained simply evaluating $dN_k/dq$, using the probability function in Eq. (\ref{pfinal}), and looking for the global maximum $q^{*}$. The result is given by
\begin{equation}
q^{*}=\frac{k(j-1)}{k\,j-1}.
\end{equation} 
Since the number of generated random trees $N_k$ is minimum, and consequently the statistical error maximum, for trees with large number of branches, the optimal value of $q$ can be obtained considering in particular the limiting case $lim_{k\to\infty} q^{*}$, yielding $q_{opt}=(j-1)/j$. Note that this coincides precisely
with the minimum value from the range of allowed values obtained in Eq. (\ref{q}).

\begin{figure}[ht] 
\begin{center}
\hbox{
\includegraphics[width=0.41\textwidth,angle=-90]{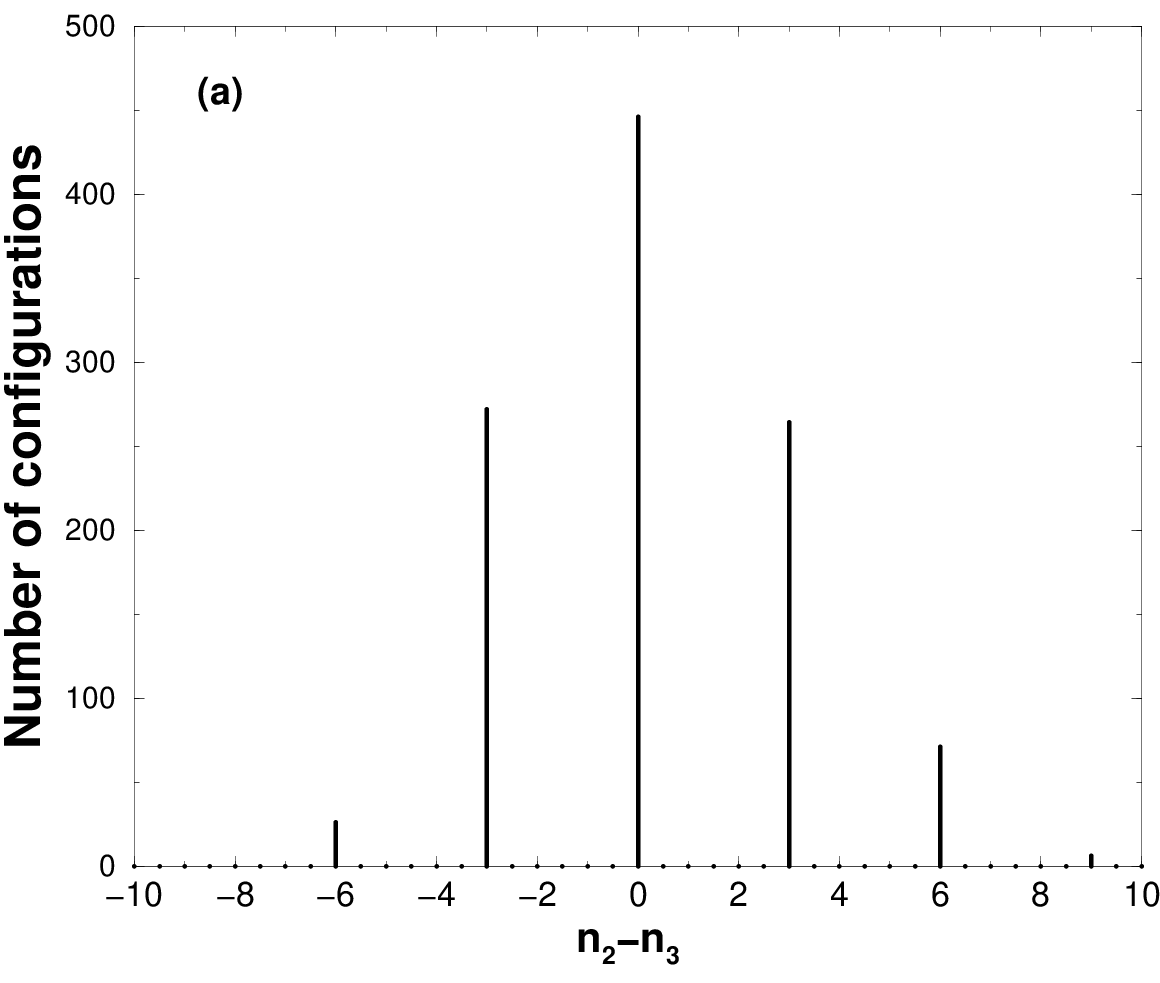}
\includegraphics[width=0.41\textwidth,angle=-90]{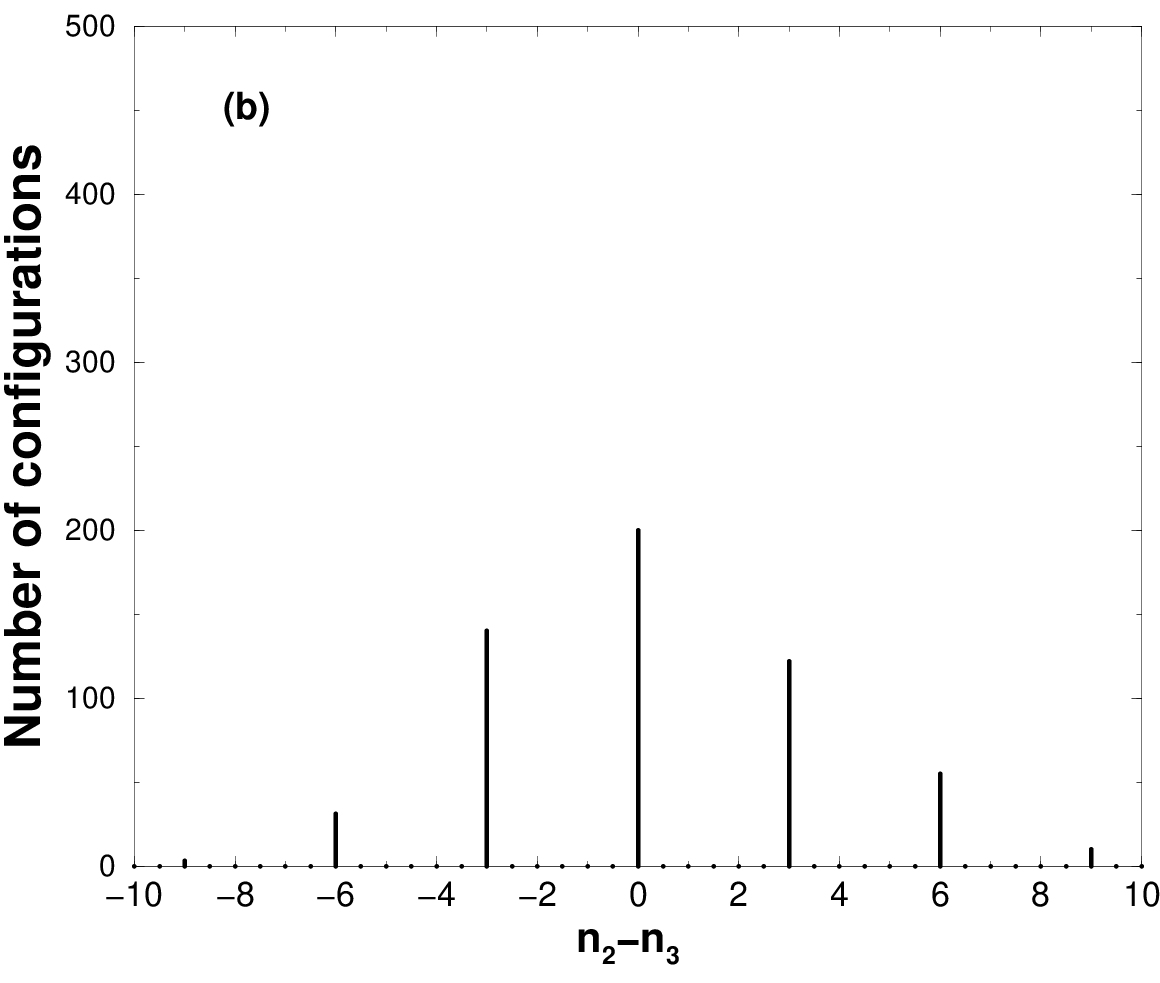}}
\caption{Number of configurations of generated random trees as function
of $n_2-n_3$, being $n_2,n_3$ the number of vertices with two and three
children, respectively. This has been done for trees with (a) $k=13$, and 
(b) $k=19$. Parameters are $N=10^5$, and $m=3$} 
\label{numberconfigurations} 
\end{center}
\end{figure}

This result can be generalized further for an arbitrary nonlinear function 
$f(u)=q+\frac{(1-q)}{m-1}\sum_{i=2}^{m} u^{i}$. Recall that for such a function
the random trees obtained may be composed in general of $m-1$ different type of vertices, each one possessing from $2$ to $m$ children. A crucial hint is to realize that the most probable configuration of any generated random tree with arbitrary number of branches $k$ occurs when the number of different type of vertices are identical, or in other words the number of children in the random tree are uniformly distributed. 
To illustrate through numerical simulations the observation above,  
$N$ random trees are generated for a nonlinear function with $m=3$, giving rise
therefore to trees composed of vertices with two and three children.  For each tree the number of vertices with two children $n_2$, and with three children $n_3$, were recorded, and computed the difference between them. In Fig. \ref{numberconfigurations} a histogram showing the number of configurations obtained as a function of $n_2-n_3$ is shown, and this has been done for two different values of the number of branches $k$. Note, as expected, that the maximum number of configurations corresponds to generated trees composed of the same number of vertices with two, and three children, that is, when $n_2-n_3=0$.

For the general function $f(u)$ with $m-1$ nonlinear terms the most probable
configuration of a random tree with $k$ branches satisfies 
$n_2=n_3=\ldots n_m=n$, being $n_j$ the number of vertices with $j$ children. Then, the global number of splitting events $Ne$ for such a configuration is given by $Ne=(m-1)n$. Clearly the number of branches $k$ of such a tree is related  with the number of splitting events $Ne$, and this relation can be obtained 
readily, yielding for $k$
\begin{equation}
k=1+\sum_{i=1}^{m-1}(m-i)n=n\frac{m-1}{2}m+1.
\end{equation} 
Given $Ne=(m-1)n$, it holds that $k=\frac{Ne}{2}m+1$. For a given finite sample size N, the number of random trees $N_k$ with $k$ branches should be maximal, and consequently the statistical error minimal, provided that the value of $q$ is chosen such as $(dN_k/dq)|_{q^{*}}=0$. This yields,
\begin{equation}
q^{*}=\frac{k\,m}{k(m+2)-2}.
\end{equation}  
Similarly to the simple case analyzed above, the optimal value of $q$ can be obtained considering in particular the limiting case $lim_{k\to\infty} q^{*}$, yielding now $q_{opt}=m/(m+2)$. Note that such a value
coincides again with the smallest value from the range of allowed values in Eq. 
(\ref{qgen}).

Apart from the errors discussed above,  a new source of error now should be taken into account. This consists of the numerical error made in approximating divergent series by a Pade approximant, since for the class of problems considered in this paper, both strategies proposed require dealing with series that in general may be divergent.
Since finding theoretical estimates of such an error for any given problem
may be a
formidable task, our goal in this section consists merely to gain some insight of such an error, illustrating how well Pade approximation actually works by analyzing a few relevant test problems. More specifically, our aim is twofold. 

On one hand we show that the statistical error made in
evaluating the partial contributions to the solution in Eq. (\ref{eq:probexpB2}) by Monte Carlo propagates to the coefficients of the Pade approximant. However, when
the statistical error is sufficiently small, the coefficients of the Pade approximant can be obtained within a reasonable accuracy, being therefore the Pade approximant rather robust at least for the examples here considered. 

On the other hand, the convergence of the Pade approximant to the solution is 
analyzed for such examples. 
Since for computational purpose the expansion generated with both strategies must be truncated, it becomes essential to determine whether the Pade approximant 
converges rapidly to the solution for the finite number of terms involved,
 or rather it is required to increase further
the order of the approximation by considering more terms in the expansion. In \cite{Angel2},
a pruning technique of the full random tree was proposed, which in practice amounts
to keep only few trees possessing a certain number of branches. This is because it was observed  that
truncating the expansion in (\ref{nonlinear_prev}) up to only a certain number 
of branches, might not affect appreciably the result, since the partial contributions to the global solution decay very rapidly as the number of branches increases. However, for the class of problems discussed here, the expansion in Eq.(\ref{eq:probexpB2}) may give rise
to a divergent series with coefficients, associated to the partial contributions to the solution, growing as the number of branches increases. Truncating the expansion, or equivalently pruning the trees, 
might be applied, but because it turns out to be uncontrollable, an special care should be taken. In particular the
effect of including more coefficients in the expansion, increasing further the order of the Pade approximation, will be analyzed for the test examples considered below. 
Recall that in general the convergence of the Pade approximant can be affected by artificial poles present in the denominator of the approximant, but not being
own by the function to be approximated, see e.g. \cite{Bender}. Therefore, to 
assess properly the validity of our findings,  it becomes essential to compute
the Pade approximant for different number of coefficients.

Concerning the apparent robustness of the Pade approximant against the statistical error affecting
the coefficients of the power series in (\ref{eq:probexpB2}), a main reason
could be that the solution of the test examples seems to be apparently locally Lipschitz. Thus, the error made in computing the coefficients of the Pade approximant should be bounded. In fact, in \cite{Wuytack} it has been proved the following related theorem
\begin{equation}
\parallel P_f-P_{f'}\parallel\leq K \parallel c-c'\parallel, 
\end{equation}
provided that $\parallel c-c'\parallel\leq d$. Here $P_f$, and $P_{f'}$  are the Pade approximants of order $(m,n)$ in $[a,b]$ of a given power series $f$ and $f'$ with coefficients $c_j$, and $c'_j$ respectively, being $\parallel c\parallel=max_{i\leq i\leq n+m} |c_i|$, $f$ locally Lipschitz, and $K$ and $d$ constants depending only on $c_i$ and $[a,b]$.

Moreover, it is worth to observe that both errors described above may be alleviated in any case by increasing conveniently the sample size $N$, and considering more coefficients in the expansion in order to compute the Pade approximant.

In the following, we present several test examples concerning one-dimensional initial value parabolic problems to illustrate what it was described before. 
The solution was computed 
probabilistically at the points $(x,t)$, where $x=0$, and $t$ several values
chosen to be distributed  between $0$ and $T$. In absence of an analytical solution the results were compared with the solution obtained upon applying an implicit finite difference scheme with a very fine mesh, and solving the ensuing algebra linear problem, characterized by a banded matrix, with LAPACK.

\noindent {\bf Example 1}. An IV parabolic problem with a purely quadratic 
negative nonlinear term. Consider the problem

\begin{eqnarray}
     u_t = u_{xx} - u^2,  \quad x\in {\bf R},\nonumber\\
     u(x,0) = \frac{e^{-\frac{x^2}{4 (t+1)}}}{\sqrt{4 \pi (t+1)}}
                                                            \label{ex1}
\end{eqnarray}


%

\begin{figure}[ht] 
\begin{center}
\hbox{
\includegraphics[width=0.4\textwidth,angle=-90]{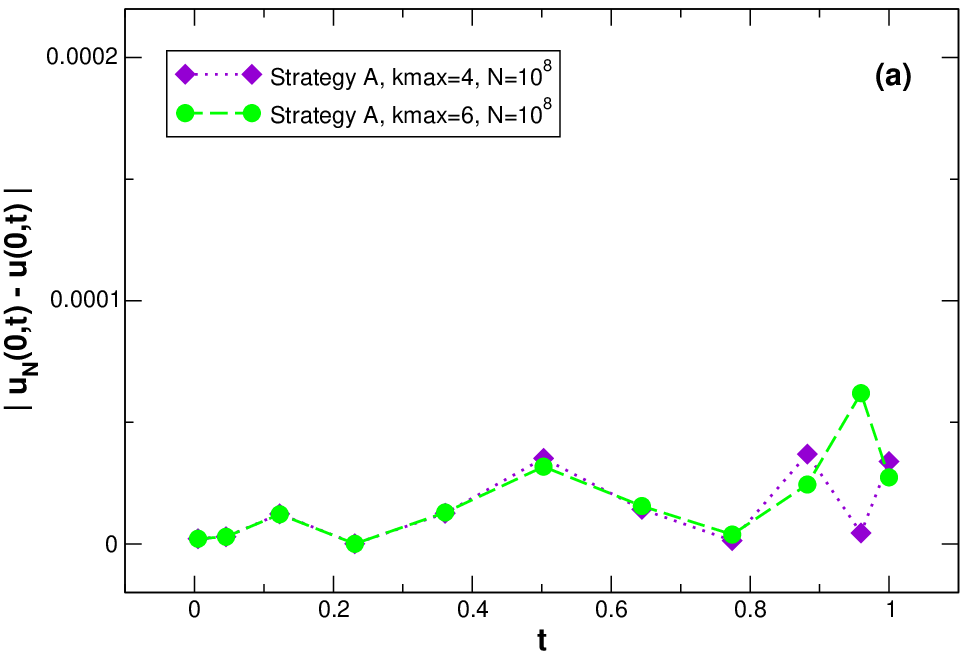}
\includegraphics[width=0.4\textwidth,angle=-90]{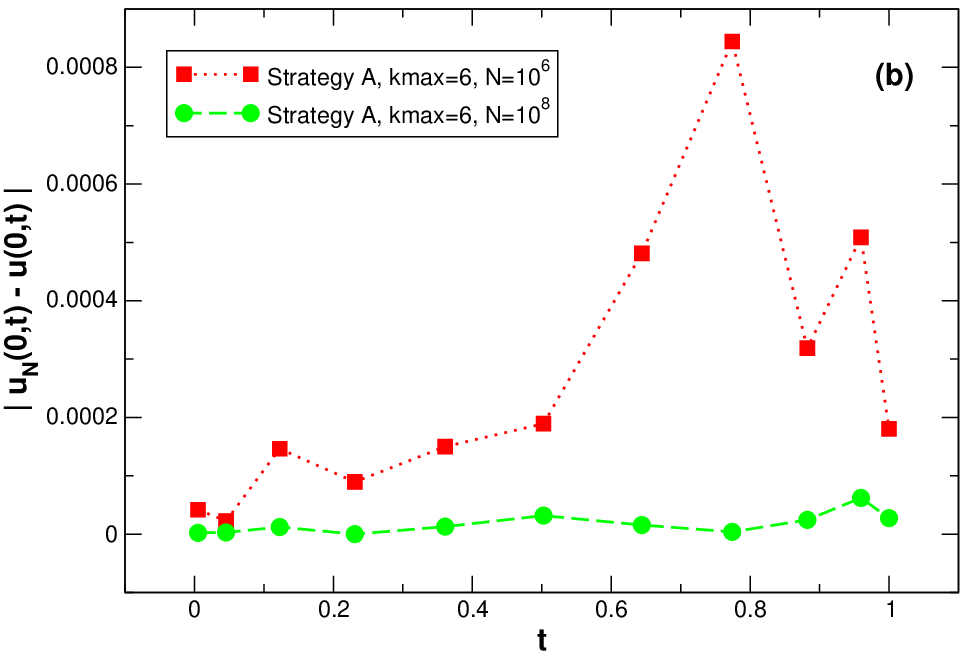}
}
\caption{Numerical error made when solving Example 1 with the strategy A.This has been done for two different values of the number of coefficients in the corresponding expansion, denoted by $kmax$.  The number of realizations $N$ was kept fixed in (a), while in (b) $kmax$ was kept fixed.}
\label{fig:ex1_errors_sA}
\end{center} 
\end{figure} 

\begin{figure}[ht] 
\begin{center}
\hbox{
\includegraphics[width=0.4\textwidth,angle=-90]{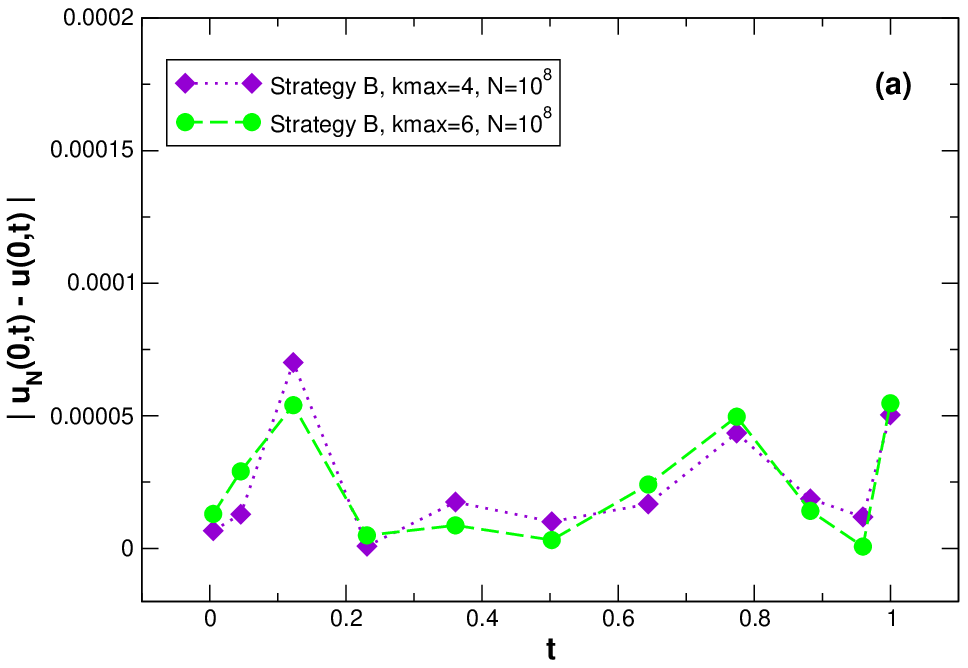}
\includegraphics[width=0.4\textwidth,angle=-90]{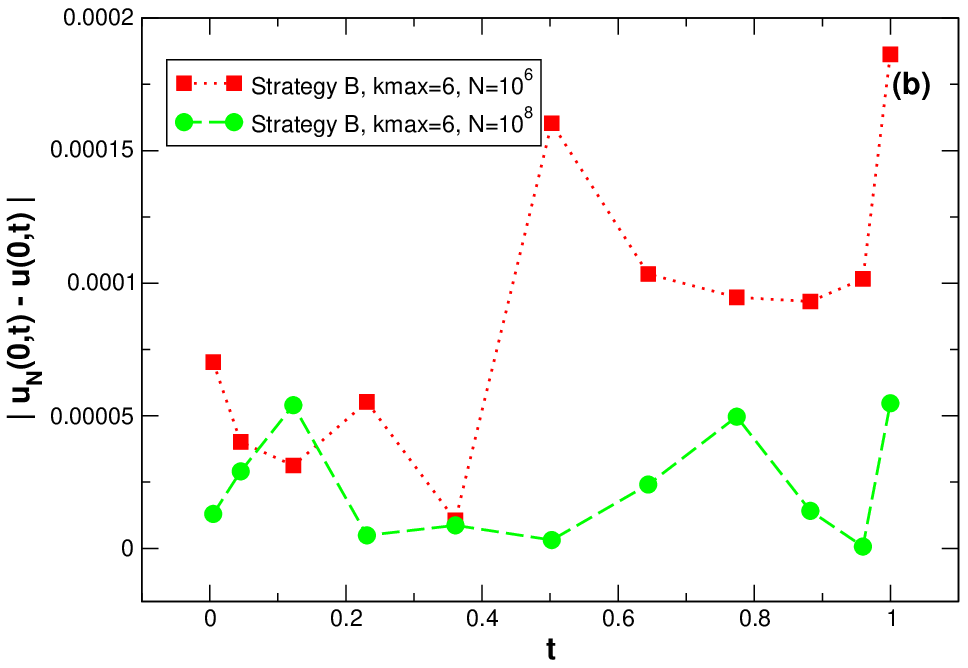}
}
\caption{Numerical error made when solving Example 1 with the strategy B. 
Identical values as in Fig. \ref{fig:ex1_errors_sA} have been used to obtain 
the results plotted in (a) and (b).}
\label{fig:ex1_errors_sB}
\end{center}
\end{figure}

The numerical error made when solving probabilistically Example 1 at a few points with  $x=0$, and $t \in [0,1]$, using both strategies, A and B, are depicted in Fig.~\ref{fig:ex1_errors_sA} and \ref{fig:ex1_errors_sB}, respectively. Note 
that for both strategies, truncating the expansion to only four coefficients, that is pruning the trees to $kmax=4$ branches, is already close to convergence for any purpose. 
Although for this example the number of coefficients to be included in the expansion could be any number above $kmax=4$, it becomes clear that choosing a larger number rather than improving accuracy, it acts reversely degrading them. In fact larger number of coefficients corresponds to contributions to the solution coming
from random trees with large number of branches, and as it was explained above such contributions are affected by larger statistical error. Moreover, it turns out to be disadvantageous as well under
a computational point of view, since generating trees with large number of branches have been proved to be rather inefficient.

Finally, note that keeping fixed the number of coefficients, and increasing the sample size, $N$, reduces accordingly the statistical error as expected.

Similar results are shown for the strategy B in Fig.~\ref{fig:ex1_errors_sB}, and therefore identical conclusions hold for this case.

\noindent {\bf Example 2}. An IV problem with a negative initial condition, 
$u(x,0)< 0$. Consider the problem

\begin{eqnarray}
     u_t = u_{xx} + u^2,  \quad x\in {\bf R},\nonumber\\
     u(x,0) =  - \, 6 \, \frac{e^{-\frac{x^2}{4 (t+1)}}}{\sqrt{4 \pi (t+1)}}
                                                            \label{ex2}
\end{eqnarray}

\begin{figure}[ht] 
\begin{center}
\hbox{
\includegraphics[width=0.4\textwidth,angle=-90]{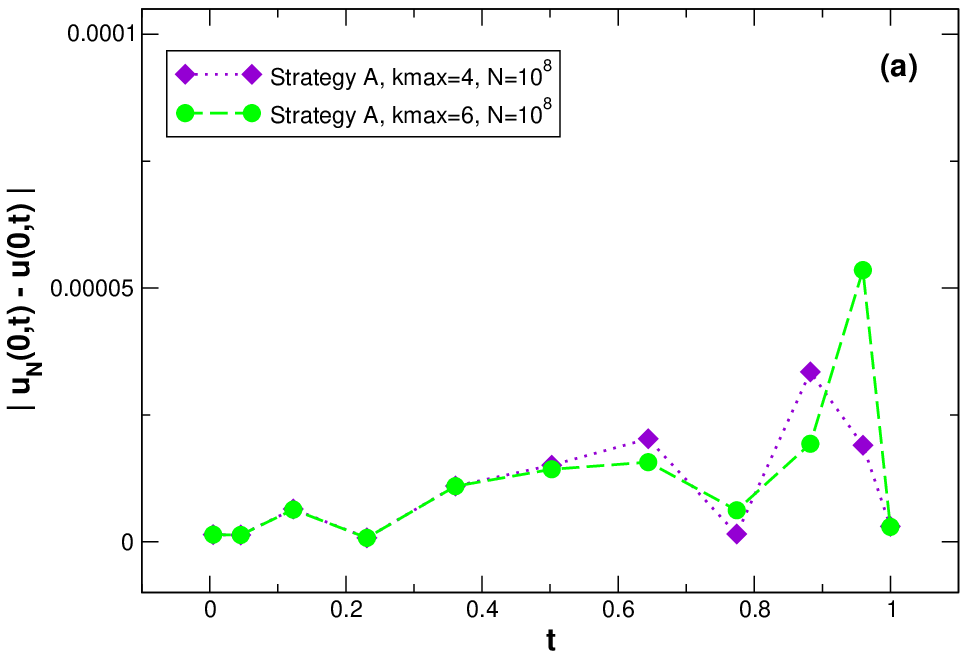}
\includegraphics[width=0.4\textwidth,angle=-90]{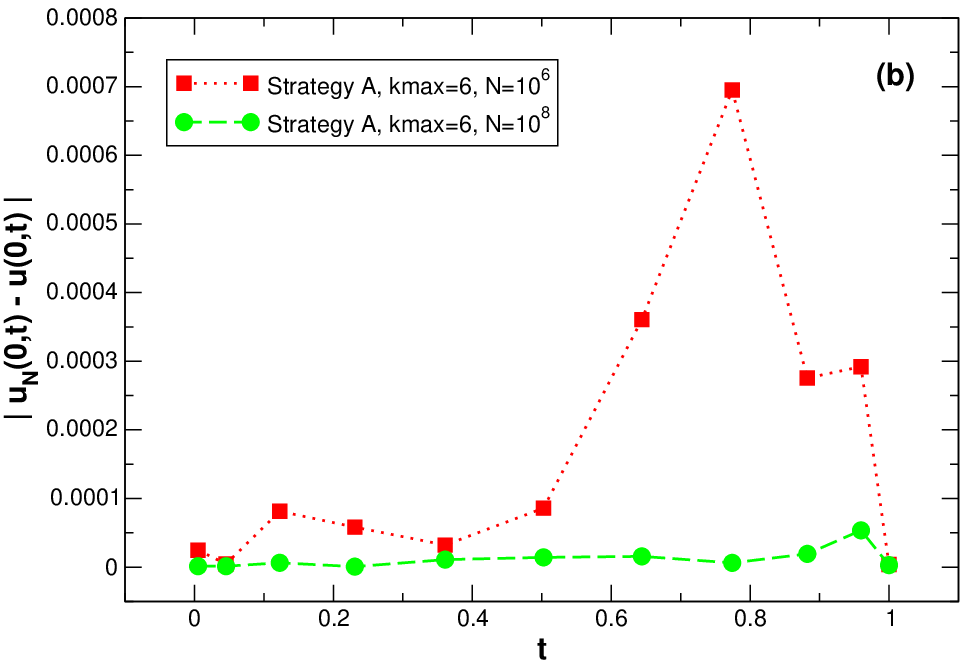}
}
\caption{Numerical error made when solving Example 2 with the strategy A. This has been done for two different number of coefficients in the corresponding expansion, denoted by $kmax$. The number of realizations $N$ was kept fixed in (a), while in (b) it was kept fixed $kmax$.}
\label{fig:ex2_errors_sA}
\end{center} 
\end{figure}

\begin{figure}[ht] 
\begin{center}
\hbox{
\includegraphics[width=0.4\textwidth,angle=-90]{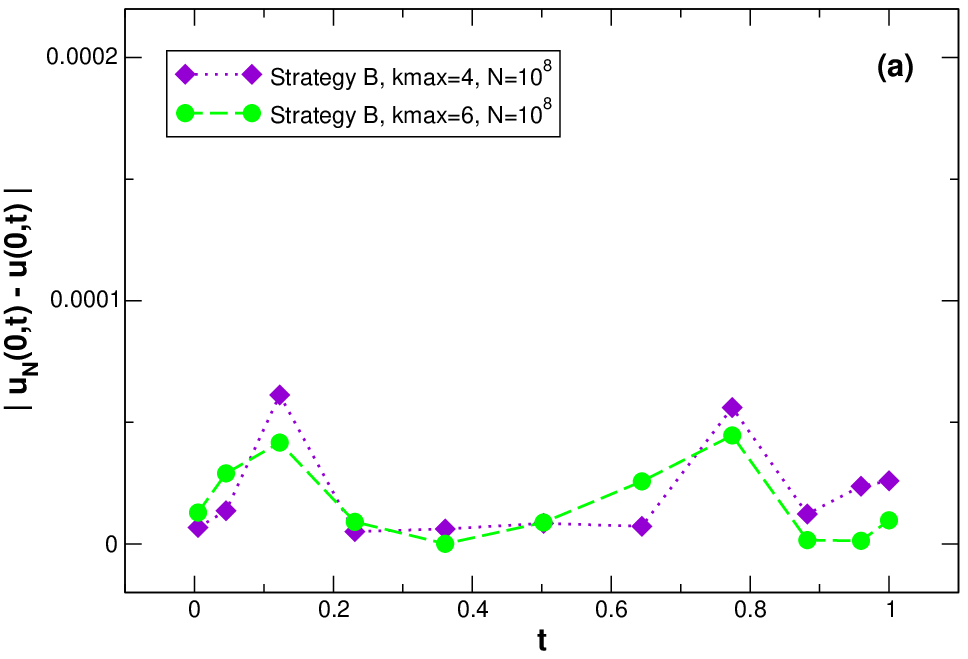}
\includegraphics[width=0.4\textwidth,angle=-90]{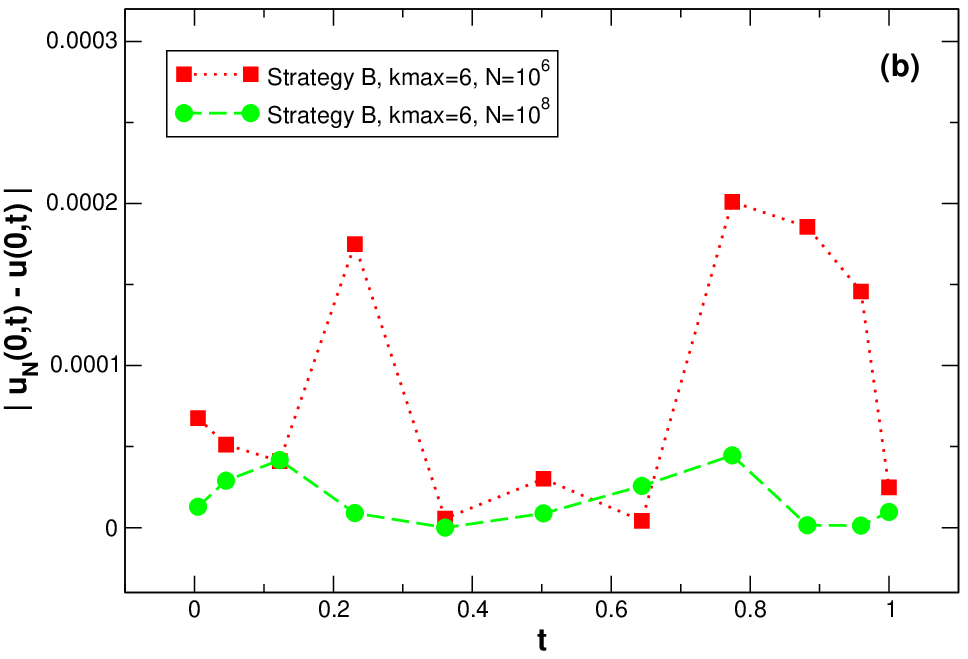}
}
\caption{Numerical error made when solving Example 2 with the strategy B. 
Identical values as in Fig. \ref{fig:ex2_errors_sA} have been used to obtain 
the results plotted in (a) and (b).}
\label{fig:ex2_errors_sB}
\end{center} 
\end{figure}

Note that the initial condition is now defined negative, and greater than $1$ in absolute value. Results are depicted in Fig.~\ref{fig:ex2_errors_sA} and \ref{fig:ex2_errors_sB}, corresponding to strategies A and B, respectively. As in Example 1, similar conclusions can be reached.

\noindent {\bf Example 3}. An IV problem with two nonlinear terms. Consider the more general problem

\begin{eqnarray}
     u_t = u_{xx} - (1 + a) u^2 - u^3,  \quad x\in {\bf R},\nonumber\\
     u(x,0) =  \frac{1}{1 + e^{- \frac{x + \sqrt{2} \left(0.5 - a\right) t}{\sqrt{2}} } },
                                                            \label{ex3}
\end{eqnarray}

where the parameter $a$ has been chosen arbitrarily to be $0.25$. 

%

\begin{figure}[ht] 
\begin{center}
\hbox{
\includegraphics[width=0.4\textwidth,angle=-90]{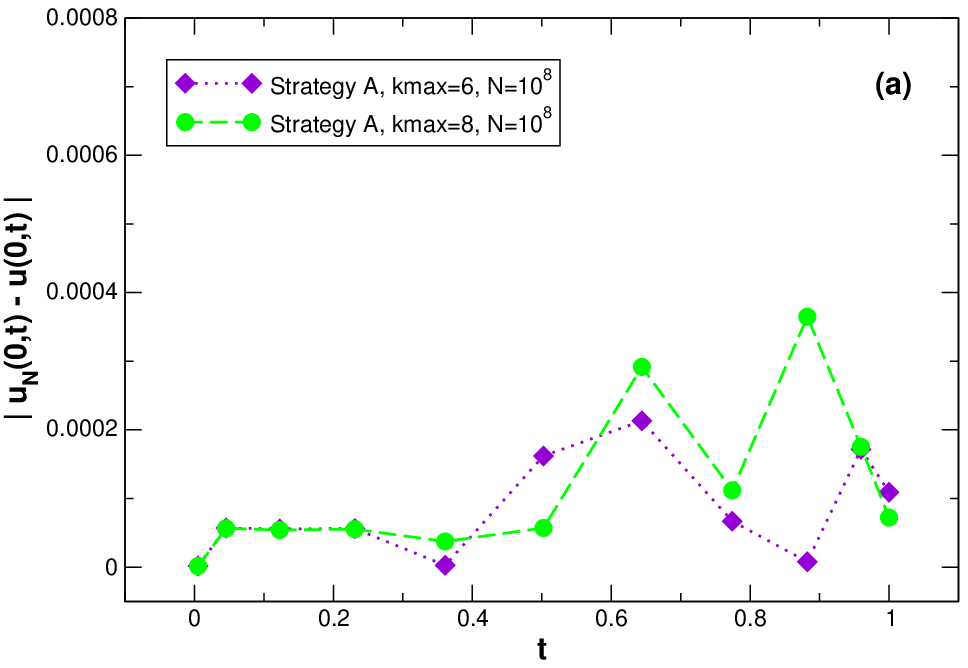}
\includegraphics[width=0.4\textwidth,angle=-90]{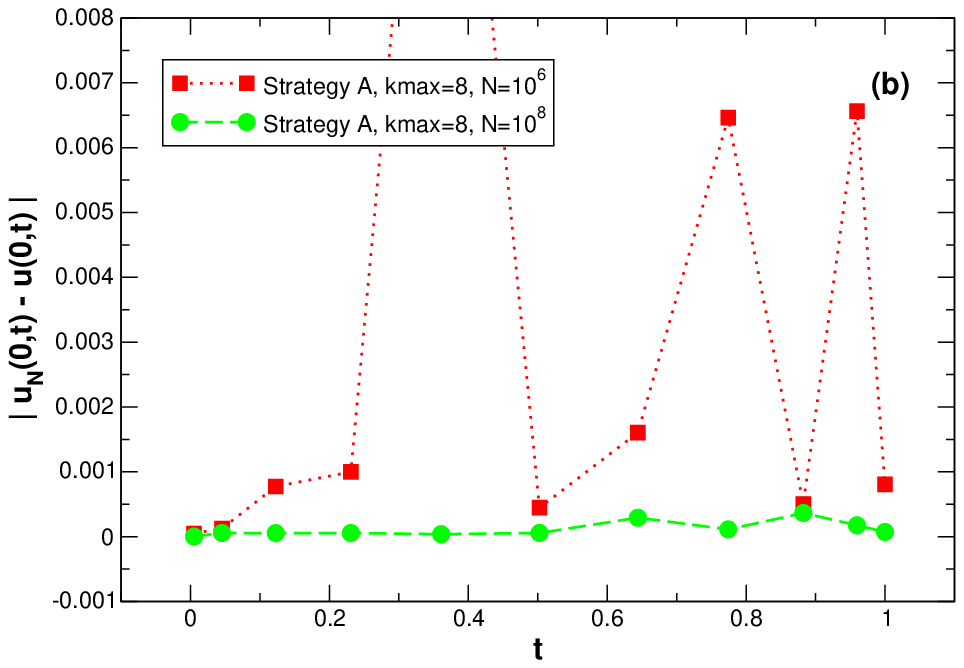}
}
\caption{Numerical error made when solving Example 3 with the strategy A. This has been done for two different number of coefficients in the corresponding expansion, denoted by $kmax$. The number of realizations $N$ was kept fixed in (a), while in (b) it was kept fixed $kmax$.}
\label{fig:ex3_errors_sA}
\end{center} 
\end{figure} 

\begin{figure}[ht] 
\begin{center}
\hbox{
\includegraphics[width=0.4\textwidth,angle=-90]{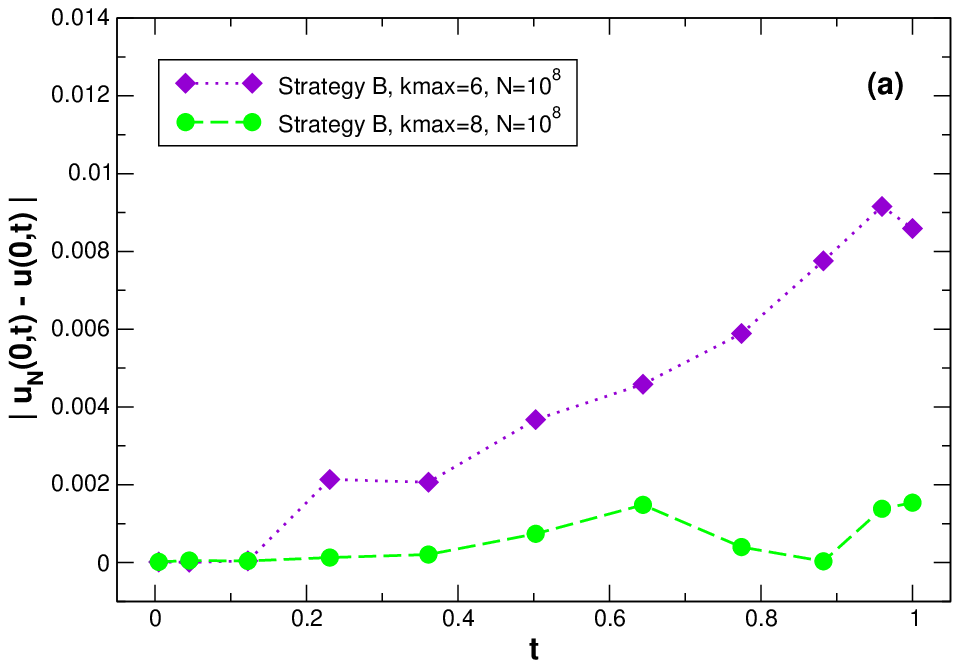}
\includegraphics[width=0.4\textwidth,angle=-90]{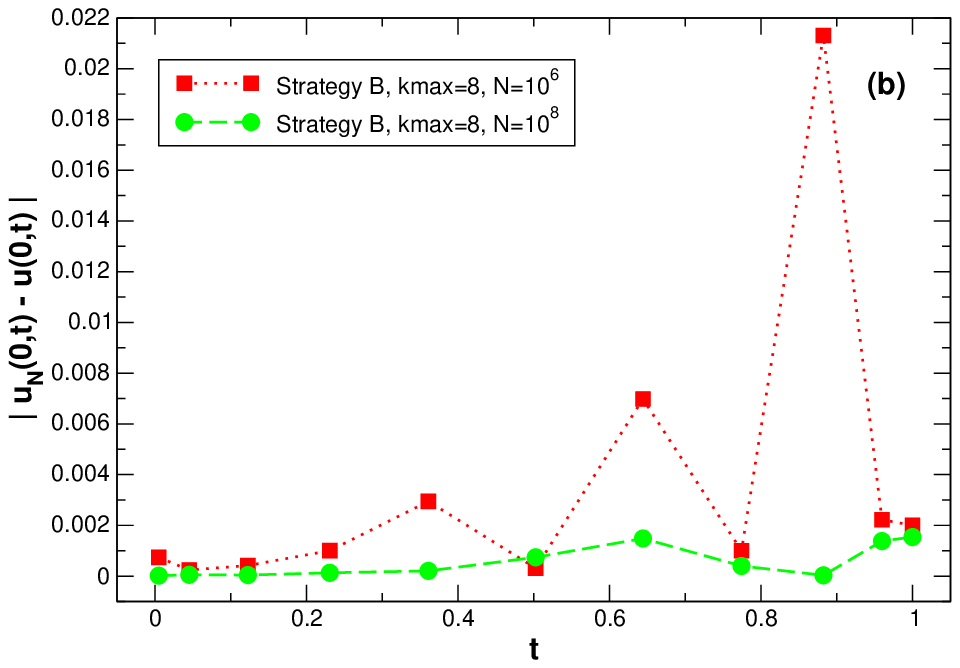}
}
\caption{Numerical error made when solving Example 3 with the strategy B. 
Identical values as in Fig. \ref{fig:ex3_errors_sA} have been used to obtain 
the results plotted in (a) and (b).}
\label{fig:ex3_errors_sB}
\end{center} 
\end{figure} 

Clearly this consists of a more involved example compared with the previous cases, since now the nonlinear function is composed of two different terms. Moreover, the coefficients multiplying both terms appear to be negative, and one of them even greater than $1$. Obviously, the joint effect of both terms gives rise to a more complex solution, suggesting the need of considering a larger number of coefficients in the expansion in order to reach convergence for the Pade approximant. This is indeed what it is observed in Fig.~\ref{fig:ex3_errors_sA}(a) and \ref{fig:ex3_errors_sB}(a). In particular for this example, it can be seen that the strategy B seems to require more coefficients than the strategy A.

Again as in the examples above, in Figs. ~\ref{fig:ex3_errors_sA}(b) and \ref{fig:ex3_errors_sB}(b) it can be observed that increasing the sample size, $N$, for both strategies reduces largely the statistical error, and in turn improves the convergence of the Pade approximant to the solution.

To conclude, both strategies showed similar performance in all test examples analyzed, however the strategy B turns out to be advantageous in any case, because 
when implemented in practice, the computational time increases linearly with the final time, growing unboundedly rather for the strategy A. 


\section{Numerical examples}
\label{examples}

The probabilistic representation described in Sec. \ref{agpr} can be hardly used for solving efficiently semilinear parabolic problems
in a whole domain, due to the high computational cost of evaluating the solution at single points. However, such a 
representation can be combined successfully with a classical domain decomposition method, as it was proposed in \cite{Acebron1,Acebron2}. The method was called probabilistic domain decomposition (PDD for short), and consists of a hybrid algorithm which requires
generating only few interfacial values along given, possibly 
artificial interfaces inside the domain, then obtaining approximate 
values upon interpolation on such interfaces. 
Such values are used as boundary data to split the original problem into a
number of fully decoupled sub-problems. 
The main advantage of this method is that the corresponding codes 
are especially suited for {\it massively parallel computing} \cite{Arbenz}. 
In fact, being the solution obtained probabilistically through an expected 
value over a given finite sample whose elements are independent from 
each other, and then after the domain decomposition the corresponding sub-problems fully decoupled, the implemented parallel
codes are characterized by an extremely low communication overhead among the various processors, affecting positively crucial properties such as {\it scalability} and {\it fault tolerance}. In the following, we describe briefly the main parts of the PDD algorithm, and for more details we refer the reader to \cite{Angel1},e.g.
  
\textit{Probabilistic part}. This is the first step to be carried out, and consists of computing the solution of the PDE at a few suitable points by some of the probabilistic strategies described in Sec. \ref{agpr}.

\textit{Interpolation}. Once the solution has been computed at few points on each interface, a second step consists of interpolating on such points, being used as nodal points, thus obtaining 
continuous approximations of interfacial values of the solution. For this purpose, since the examples analyzed below corresponds to two-dimensional problems  a
tensor product interpolation based on cubic spline \cite{Antia} was used. The computational cost of this part turns out to be negligible compared with the time spent in the other parts of the algorithm. 
The nodal points are uniformly distributed on each plane, and a not-a-knot condition is imposed.

\textit{Local solver}. The third and final step consists of computing the solution inside 
each subdomain, this task being assigned to different processors. 
This can be accomplished resorting to local solvers, which may use classical nume\-rical schemes, such as implicit finite differences for simple geometries or finite elements methods for more complex configurations. For the former case, subroutines based on LAPACK for solving the ensuing linear algebra problems has been chosen, since the corresponding matrices are banded. Therefore, each processor
can be devoted only to the solution of its local linear system, whose banded associated 
matrix is smaller. Concerning the memory consumption per processor, 
including an extra fill-in space, the total amount is considerably
reduced \cite{Angel2}.

\begin{figure}[ht] 
\begin{center}
\vbox{
\includegraphics[width=0.45\textwidth ,angle=0]{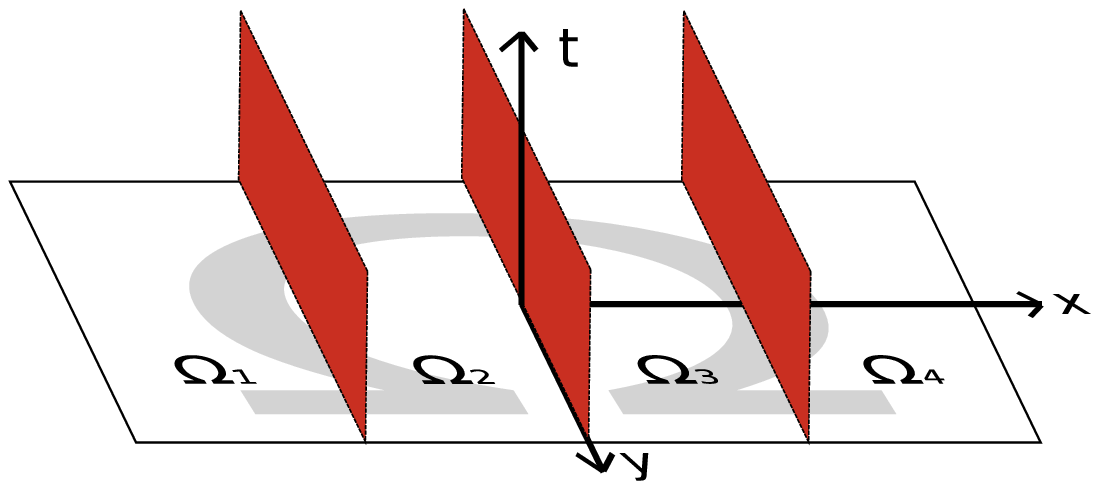}
\includegraphics[width=0.45\textwidth ,angle=0]{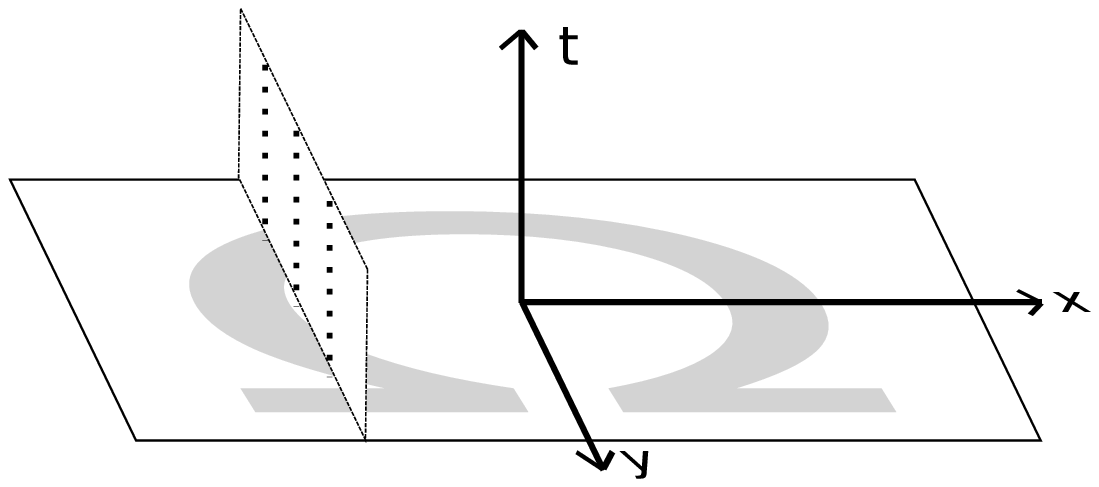}
}
\caption{A sketchy diagram illustrating the main steps of the algorithm in 2D: 
The figure on the left shows how the domain decomposition is done in 
practice. The figure on the right shows the points where the 
solution is computed probabilistically; these are used afterward
as nodal points for interpolation.}
\label{sketch_algorithm}
\end{center} 
\end{figure}

In Fig.~\ref{sketch_algorithm} we sketch a diagram, illustrating how the algorithm 
works in practice for a two-dimensional case. Here the solution is obtained probabilistically 
at a few points pertaining to some ``interfaces'' conveniently chosen inside the space-time domain 
$D:=\Omega \times [0, T]$, with $\Omega\subset \mathbb{R}^2$. Such interfaces divide 
the domain into $p$ subdomains, $\Omega_i$, from $i=1,...,p$, being assigned to different
 processors, $p_i, i=1,...,p$. The more convenient way to parallelize this part is splitting in in\-de\-pen\-dent sets of points. Since the number of points where the solution is computed is larger than the number of processors $p$, computing such a solution can be assigned as a task to different processors. This can be seen as a
 coarse-grain parallelization, and even though other finest strategies can be adopted, this one turns out to be the more convenient for the examples analyzed in this section.

  Here we present some numerical examples for 2D initial value
problems to 
illustrate the PDD
algorithm, being the probabilistic part built up with the two strategies A and B discussed in the previous section. All simulations were carried out on the Matrix supercomputer, belonging to the Inter-University Consortium for the Application of Super-Computing 
for Universities and Research (CASPUR) located in Rome (Italy), using up to 512 processors. This supercomputer consists of a Linux 
cluster based on multi-core Opteron processor nodes with Infiniband 
interconnection, and it was ranked in the Top500 list with a peak performance 
of 22 TFlops.

  As in \cite{Angel1,Angel2}, a comparison was made solving the same problems by some other classical numerical 
methods in order to asses the performance of both methods. For the space-time domain as well as for the subdomains in our decomposition, we used the Crank-Nicolson 
(implicit) finite difference method. On the various decoupled subdomains obtained by the PDD algorithm we 
used LAPACK for solving the ensuing linear algebra problem, while the full domain solution was computed by ScaLAPACK. This widely used and freely available numerical package has been considered extremely 
efficient for the parallel solution of banded linear systems. For more details concerning the computational 
cost of both methods, LAPACK and ScaLAPACK, see \cite{Angel2}, e.g.

\noindent {\bf Example 4}. An IV problem with two nonlinear terms. Consider the problem

\begin{eqnarray}
     u_t &=& u_{xx} + u_{yy} - (1 + a) u^2 - u^3,  \quad (x,y)\in {\bf R}^2,\nonumber\\
     && u(x, y, 0) = -2 \cos^2(\frac{\pi\,x } {2 A_x}) 
\cos^2(\frac{\pi\, y}{2 A_y}).
\end{eqnarray}

where $a=0.25, A_x = 10, A_y = 40$. The space and
time discretization step has been chosen to be $\Delta x = \Delta y = 0.25, \Delta t=10^{-3}$, and the solution was computed for a final time $T = 0.5$.

Note that the unbounded domain should be truncated conveniently to a boun\-ded
domain in order to be able to solve numerically the problem using a finite difference scheme. This requires introducing some artificial boundary conditions to confine the computational domain. Since the problem is formulated as a pure initial value problem, the artificial boundary conditions should be prescribed in such a way no additional data are imposed on such boundaries. In practice, this can be done readily imposing Dirichlet boundary conditions on the artificial boundaries, such that the boundary conditions are automatically satisfied by the solution of the problem. However, being the solution of the problem unknown, one should 
resort to several type of approximations of the solution to be used as boundary conditions.
For the problem above, the solution of the problem is assumed to decay sufficiently fast to infinity, and being the computational domain chosen to be large enough $\Omega\in [-L_x,L_x]\times[-L_y,L_y]$, with $L_x=40$ and $L_y=160$, a zero Dirichlet boundary condition can be properly imposed at $x=\pm L_x$,$y=\pm L_y$. 

When a probabilistic representation is available,  such a representation were used as well to obtain much more accurate approximations for the artificial boundary conditions, since it allows to obtain the solution at any single point arbitrarily chosen. This is remarkable feature of the probabilistic representation, not owned by any other numerical method.

\begin{figure}[ht] 
\centering 
\includegraphics[width=0.9\textwidth,angle=0]{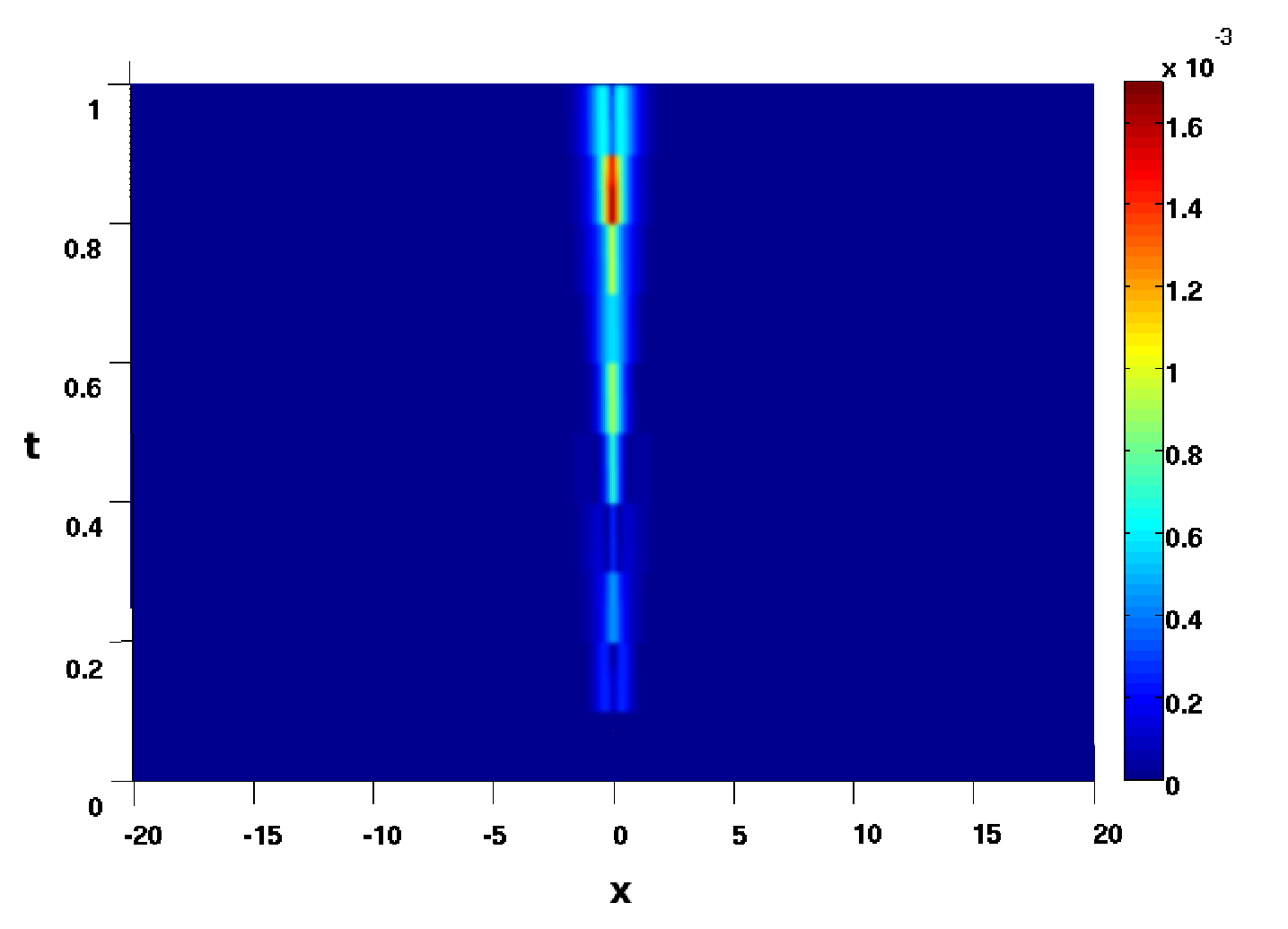}
\caption{Pointwise numerical error made when solving Example 4 with the strategy A.}
\label{fig:pointwise_error_SA}
\end{figure}

\begin{figure}[ht] 
\centering
\includegraphics[width=0.94\textwidth,angle=0]{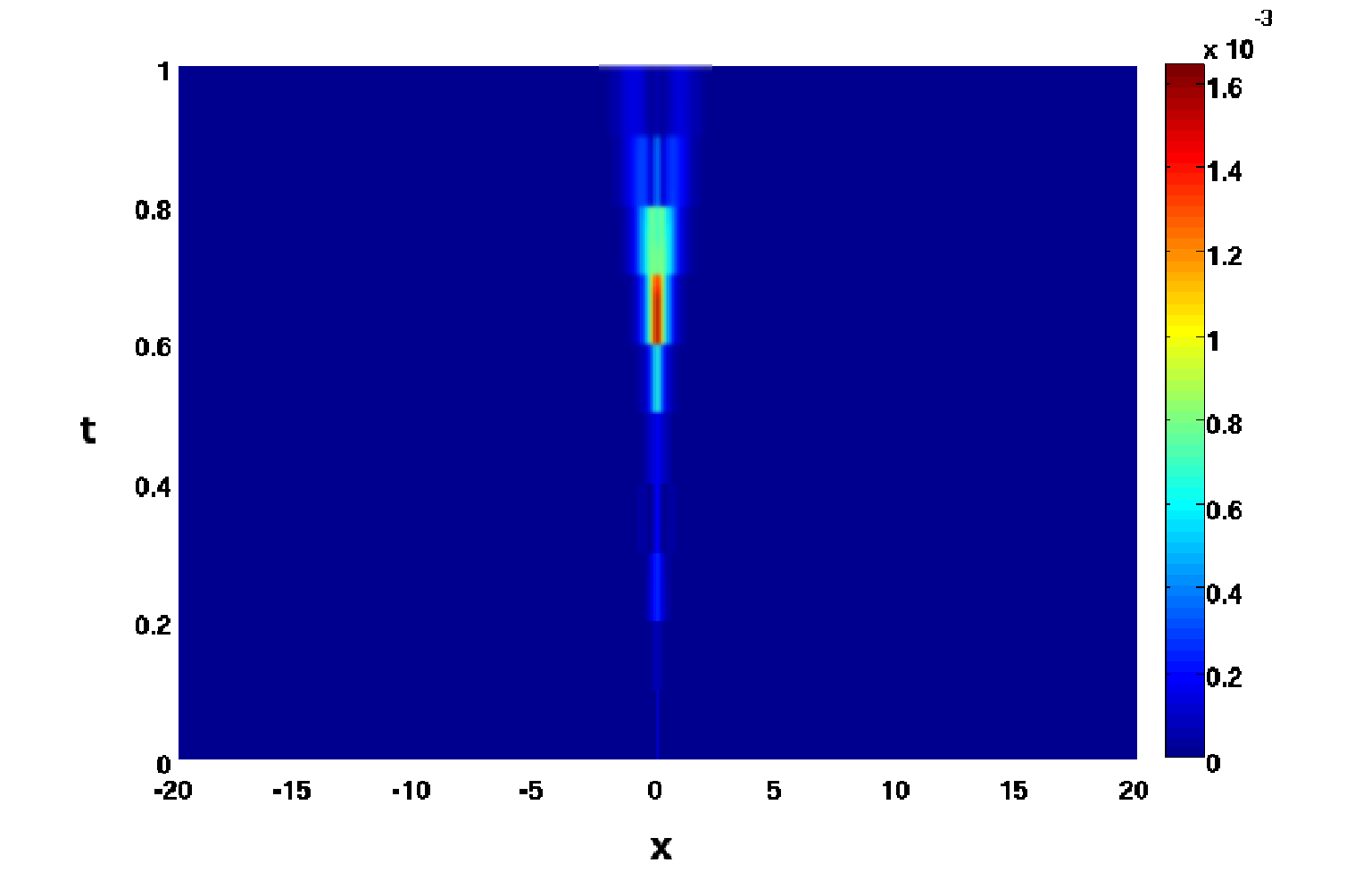}
\caption{Pointwise numerical error made when solving Example 4 with the strategy B.}
\label{fig:pointwise_error_SB}
\end{figure}

In Fig.~\ref{fig:pointwise_error_SA} and Fig.~\ref{fig:pointwise_error_SB} the pointwise numerical error around the origin made with the strategy A and strategy B, respectively, is shown. Here the value of the
parameters were kept fixed to $\Delta x = \Delta y = 10^{-2}$, $\Delta t = 10^{-3}$, $T=1$. For clarity, only the maximum absolute value of the error obtained in the $y-$axis is plotted. 

\begin{table}[htbp]   
\caption{Example 4: $T_{PDD}$, $T_{ScaLAPACK}$ denotes the computational time spent in seconds by PDD with the strategy A, and ScaLAPACK, respectively. 
$T_{MC}$, and $T_{INT}$ correspond to the time spent by the probabilistic and the interpolation part, respectively; $Memory$ denotes the total memory consumption.}
\begin{center}
\begin{tabular}%
     {| @{\hspace{5mm}} c @{\hspace{5mm}}|| @{\hspace{5mm}} c @{\hspace{5mm}} c @{\hspace{10mm}} r@{.}l@{ GBs} @{\hspace{5mm}} c || c |}
\hline
    \textbf{\em Procs.}
  & \multicolumn{1}{l}{$T_{MC}$}
  & \multicolumn{1}{l}{$T_{INT}$}
  & \multicolumn{2}{l}{\textit{Memory}}
  & $\mathbf{T_{PDD}}$
  & $\mathbf{T_{ScaLAPACK}}$  \\ \hline \hline
\textbf{128}  &  902\textquotedblright & $<$1\textquotedblright &  0&86 & 5572\textquotedblright & 29881\textquotedblright \\
\textbf{256}  &  998\textquotedblright & $<$1\textquotedblright &  0&19 & 2086\textquotedblright & 23953\textquotedblright \\
\textbf{512}  &  1018\textquotedblright & $<$1\textquotedblright &  0&06 & 1327\textquotedblright & 23334\textquotedblright \\
\hline
\end{tabular}
\label{TableA}
\end{center}

\end{table}

\begin{table}[htbp]   
\caption{Example 4: $T_{PDD}$, $T_{ScaLAPACK}$ denotes the computational time spent in seconds by PDD with the strategy B, and ScaLAPACK, respectively. 
$T_{MC}$, and $T_{INT}$ correspond to the time spent by the probabilistic and the interpolation part, respectively; $Memory$ denotes the total memory consumption.}
\begin{center}
\begin{tabular}%
     {| @{\hspace{5mm}} c @{\hspace{5mm}}|| @{\hspace{5mm}} c @{\hspace{5mm}} c @{\hspace{10mm}} r@{.}l@{ GBs} @{\hspace{5mm}} c || c |}
\hline
    \textbf{\em Procs.}
  & \multicolumn{1}{l}{$T_{MC}$}
  & \multicolumn{1}{l}{$T_{INT}$}
  & \multicolumn{2}{l}{\textit{Memory}}
  & $\mathbf{T_{PDD}}$
  & $\mathbf{T_{ScaLAPACK}}$  \\ \hline \hline
\textbf{128}  &  736\textquotedblright & $<$1\textquotedblright &  0&86 & 5413\textquotedblright & 29881\textquotedblright \\
\textbf{256}  &  824\textquotedblright & $<$1\textquotedblright &  0&19 & 1917\textquotedblright & 23953\textquotedblright \\
\textbf{512}  &  837\textquotedblright & $<$1\textquotedblright &  0&06 & 1152\textquotedblright & 23334\textquotedblright \\
\hline
\end{tabular}
\label{TableB}
\end{center}
\end{table}

In Table~\ref{TableA} and Table~\ref{TableB}, the computational times obtained when solving the example 4 using the PDD algorithm with the strategy~A and the strategy~B, respectively, are shown. The partial computational times spent by the probabilistic part and the interpolation part of the algorithm, as well as the computational time spent by ScaLAPACK, have also been displayed. The two methods were
compared correspondingly to the same maximum error, $10^{-3}$.

It is worth to observe that the computational times obtained with the strategy B are significantly smaller than those obtained with the strategy A. This can be explained in view of the less computational cost of the probabilistic representation based on the strategy B, as it was already theoretically shown in the previous section.

\noindent {\bf Example 5}. An IV problem with a single nonlinear term and variable coefficients. Consider the problem

\begin{eqnarray}
     u_t &=& u_{xx} + u_{yy} - e^{-x^2} \frac{1}{t+0.1} u^2,  \quad (x,y)\in {\bf R}^2,\nonumber\\
     && u(x, y, 0) = \frac{e^{-\frac{x^2}{4}}}{\sqrt{4 \pi}}. 
\end{eqnarray}

The space and time steps are $\Delta x = \Delta y = 0.25, \Delta t=10^{-3}$, and the solution was computed at the final time $T = 0.5$.

Note that in this example a variable coefficient depending on time and space, multiplying the nonlinear term $u^2$ has been considered, and it may be in general taken values larger than one.

The computational times are shown in Table~\ref{TableC} and Table~\ref{TableD}, comparing the strategy A and B, respectively, with ScaLAPACK. Note that the computational times turns out to be slightly smaller than those obtained in the previous example, and this is because the generated random trees now are purely binary. The strategy B wins again over the strategy A, and the same reason of
the previous example holds also for the present case.

\begin{table}[htbp]   
\caption{Example 5: $T_{PDD}$, $T_{SCALAPACK}$ denotes the computational time spent in seconds by PDD with the strategy A, and ScaLAPACK, respectively. 
$T_{MC}$, and $T_{INT}$ correspond to the time spent by the Monte Carlo and the interpolation part, respectively; $Memory$ denotes the total memory consumption.}
\begin{center}
\begin{tabular}%
     {| @{\hspace{5mm}} c @{\hspace{5mm}}|| @{\hspace{5mm}} c @{\hspace{5mm}} c @{\hspace{10mm}} r@{.}l@{ GBs} @{\hspace{5mm}} c || c |}
\hline
    \textbf{\em Procs.}
  & \multicolumn{1}{l}{$T_{MC}$}
  & \multicolumn{1}{l}{$T_{INT}$}
   & \multicolumn{2}{l}{\textit{Memory}}
  & $\mathbf{T_{PDD}}$
  & $\mathbf{T_{ScaLAPACK}}$  \\ \hline \hline
\textbf{128}  &  807\textquotedblright & $<$1\textquotedblright &  0&86 & 5461\textquotedblright & 27466\textquotedblright \\
\textbf{256}  &  889\textquotedblright & $<$1\textquotedblright &  0&19 & 1954\textquotedblright & 22070\textquotedblright \\
\textbf{512}  &  906\textquotedblright & $<$1\textquotedblright &  0&06 & 1102\textquotedblright & 21327\textquotedblright \\
\hline
\end{tabular}
\label{TableC}
\end{center}
\end{table}

\begin{table}[htbp]   
\caption{Example 5: $T_{PDD}$, $T_{ScaLAPACK}$ denotes the computational time spent in seconds by PDD with the strategy B, and ScaLAPACK, respectively. 
$T_{MC}$, and $T_{INT}$ correspond to the time spent by the probabilistic and the interpolation part, respectively; $Memory$ denotes the total memory consumption.}
\begin{center}
\begin{tabular}%
     {| @{\hspace{5mm}} c @{\hspace{5mm}}|| @{\hspace{5mm}} c @{\hspace{5mm}} c @{\hspace{10mm}} r@{.}l@{ GBs} @{\hspace{5mm}} c || c |}
\hline
    \textbf{\em Procs.}
  & \multicolumn{1}{l}{$T_{MC}$}
  & \multicolumn{1}{l}{$T_{INT}$}
  & \multicolumn{2}{l}{\textit{Memory}}
  & $\mathbf{T_{PDD}}$
  & $\mathbf{T_{ScaLAPACK}}$  \\ \hline \hline

\textbf{128}  &  511\textquotedblright & $<$1\textquotedblright &  0&86 & 5148\textquotedblright & 27466\textquotedblright \\
\textbf{256}  &  562\textquotedblright & $<$1\textquotedblright &  0&19 & 1623\textquotedblright & 22070\textquotedblright \\
\textbf{512}  &  569\textquotedblright & $<$1\textquotedblright &  0&06 & 882\textquotedblright & 21327\textquotedblright \\
\hline

\end{tabular}
\label{TableD}
\end{center}
\end{table}

\section{Summary}
The class of semilinear parabolic problems amenable to a probabilistic
solution has been expanded by introducing suitable generalized random
trees. The probabilistic computation consists of evaluating averages
on the generated random tree, which plays a role similar to that of a
random path in linear problems. The new representation allows
treatment of semilinear problems without a potential term, with
arbitrary coefficients multiplying the nonlinear term, and arbitrary
initial data, including negative definite and greater than
one. The implementation uses two different strategies, which
require computing the solution through a series where the coefficients
represent the partial contribution to the solution coming from
generated random trees with any number of branches. Since such a
series might be divergent, in general it cannot be summed simply by a
sequence of partial sums. Nevertheless, numerical experiments show
that, in many cases, the asymptotic series can be approximated quite
accurately by techniques based on the Pade
approximant. A qualitative analysis of the error done has been 
carried out, showing that for the test problems analyzed so far,
considering a few coefficients of the series suffices to obtain a
reasonable accuracy. 

 Moreover, it has been shown that the strategy
termed  B greatly reduces the computation time compared with
strategy A, which is based rather on generating random trees governed by an
exponential random time.  The new probabilistic representation has been used
successfully as a crucial element for implementing a suitable probabilistic domain decomposition 
method. In contrast to the classical
deterministic method for solving partial differential equations, the
probabilistic approach computes the solution at single points
internal to the domain, without first generating a
computational mesh and solving the full problem. The generalized
PDD method has been shown to be suited for massively parallel
computers. In fact, some numerical examples have been run that show
excellent scalability properties of the PDD algorithm in
large-scale simulations, using up to 512 processors on a high
performance supercomputer. Finally, the performance of the algorithm
has been compared with other efficient, freely available parallel
algorithms, showing a striking difference.

\section*{Acknowledgments}
 
  This work was supported by the Portuguese FCT under grant 
PTDC/EIA-CCO/098910/2008.
The authors thankfully acknowledge the computer resources, technical expertise
and assistance provided by the Rome Supercomputing center (CASPUR).

\end{document}